\documentclass[11pt]{article}
\usepackage{epsfig}
\usepackage{amssymb,amsmath,amsthm,amscd}
\usepackage{latexsym}
\usepackage{amsthm}
\usepackage{xcolor}
\usepackage{amsthm}%\usepackage[notref,notcite]{showkeys}
\usepackage{graphicx}
\usepackage{xcolor}
\newtheorem{theorem}{Theorem}[section]

\newtheorem{thm}[theorem]{Theorem}
\newtheorem{lemma}[theorem]{Lemma}
\newtheorem{lem}[theorem]{Lemma}

\newtheorem{cor}[theorem]{Corollary}

\theoremstyle{definition}

\newtheorem{definition}[theorem]{Definition}

\newtheorem{rem}[theorem]{Remark}

\newcounter{labelflag} 

\newcommand{\Label}[1]{
                       \ifnum\thelabelflag=1
                          \ifmmode
                             \makebox[0in][l]{\qquad\fbox{\rm#1}}
                          \else
                             \marginpar{\vspace{0.7\baselineskip}
                                        \hspace{-1.1\textwidth}
                                        \fbox{\rm#1}}
                          \fi
                       \fi
                       \label{#1} }

\newcommand{\be}{\begin{equation}}
\newcommand{\ee}{\end{equation}}

\begin{document}

\baselineskip=1.2\baselineskip
\pagestyle{plain}
\title{Regularity of invariant measures  of solutions for fractional stochastic partial differential equations with memory and nonlinear noise
 \footnote{The research has  been partially supported by NSFC (11971394,  12371178, 12531017, 12301234, 12471165) and Natural Science Foundation of Sichuan Province
Grant (2025ZNSFSC0068), Doctoral Innovation Fund Program of Southwest Jiaotong University (CX2025YB30,CX2026YB29), China Scholarship Council (202507000029) and Spanish Ministerio de Ciencia, Innovaci\'on y Universidades (MCIU), Agencia Estatal de Investigaci\'on (AEI) and FEDER (PID2024-156228NB-I00).} }
\author{{Ran Li$^\text{1,2}$ ~~Jiaohui Xu$^\text{3}$\footnote{Corresponding author:
jxu@us.es (J. Xu)}~~Dingshi Li$^\text{1}$~~ Tom\'{a}s Caraballo $^\text{2,4}$~~  }\\
{\small\textsl{$^\text{1}$ School of Mathematics, Southwest Jiaotong University, Chengdu, 610031, P.R. China.}}\\
{\small\textsl{$^\text{2}$ Dpto. Ecuaciones Diferenciales y An\'{a}lisis Num\'{e}rico, Facultad}}\\
{\small\textsl{de Matem\'{a}ticas, Universidad de Sevilla, c/ Tarfia s/n, 41012-Sevilla, Spain.}}\\
{\small\textsl{$^\text{3}$ Center for Nonlinear Studies, School of Mathematics, Northwest University, Xi'an, 710127,
P. R. China.}}\\
{\small\textsl{$^\text{4}$ Department of Mathematics, Wenzhou University, Wenzhou, Zhejiang Province, 325035,
P. R. China.}}
}
\date{}\maketitle
	{\noindent{\bf Abstract}}
We are concerned with stochastic fractional partial differential equations involving  memory and nonlinear multiplicative noise. Based on the well-known Dafermos method, a suitable  variable defined  by the integrated past history is introduced.
This transformation not only overcomes the non-Markovianity of the original problem which prevents us from analyzing invariant measures for the associated Markov semigroup,
but also facilitates the analysis of  the well-posedness together with the Faedo-Galerkin approximations and Minty's technique. Moreover, we prove the existence of invariant measures by the classical Krylov-Bogoliubov approach, and further establish the regularity of invariant measures by virtue of a control argument.

{\noindent{\bf Keywords:}} Nonlocal fractional stochastic differential equations; Memory terms; Nonlinear multiplicative noise;  Existence and regularity of invariant measures.

 {\noindent{\bf MSC 2010.}} Primary  37L55.  Secondary 34F05, 37L30, 60H10.

\section{Introduction}
\setcounter{equation}{0}
Let $\mathcal{O}\subset \mathbb{R}^{n}$ be a bounded open domain with smooth boundary. We study a class of semilinear nonlocal fractional evolution equations with memory  and a nonlinear multiplicative noise perturbation. More precisely,
we consider the following problem,
\begin{equation}\label{1.1}
\left\{
\begin{aligned}
&\mathrm{d}u(t)=-\kappa(-\Delta)^{\alpha}u(t)\mathrm{d}t
 -(1-\kappa)\int_{0}^{\infty}K(s)(-\Delta)^{\alpha}u(t-s)\mathrm{d}s\mathrm{d}t-\lambda u(t)\mathrm{d}t
+\varphi(u(t))\mathrm{d}t\\[0.8ex]
&+\sigma(u(t))\mathrm{d}W(t),\qquad (t,x)\in \mathbb{R}^{+}\times \mathcal{O},\\[0.8ex]
&u(t)|_{\partial\mathcal{O}}=0,\quad t\in \mathbb{R}^{+},\qquad
u(0,x)=u_{0}(x), \quad x\in\mathcal{O},\\[0.8ex]
&u(t,x)=\psi(t,x)\quad t<0,\quad x\in\mathcal{O}.
\end{aligned}
\right.
\end{equation}
$\psi$ will be introduced later. Such kind of equations arise naturally in the modeling of anomalous diffusion phenomena in
viscoelastic or memory-dependent media, in which the present state of the system is influenced
by its past history.
The operator $(-\Delta)^{\alpha}$ with $\alpha\in(0,1)$ denotes the fractional Laplacian
associated with homogeneous Dirichlet boundary conditions on $\mathcal{O}$ and describes
 nonlocal diffusion effects.
The parameter $\kappa\in(0,1)$ determines the relative contribution of  the memory term.
The kernel $K:\mathbb{R}^{+}\to\mathbb{R}^{+}$ represents the memory effect,  which is assumed to  be smooth and exponentially decaying.
The nonlinear mapping $\varphi:\mathbb{R}\to\mathbb{R}$ denotes a drift term which typically
satisfies appropriate dissipativity and growth conditions. $\lambda$ is a positive constant. The model is perturbed by a nonlinear multiplicative noise
$\sigma(u(t))\mathrm{d}W(t)$, where $\sigma$ is  a  Lipschitz nonlinear diffusion coefficient, and
$W=(W(t))_{t\ge0}$ is a $\mathcal{U}$-valued cylindrical Wiener process defined on $(\Omega,\mathcal F,\{\mathcal F_t\}_{t\geq0},\mathbb P)$.

Stochastic partial differential equations arise in a wide range of physical applications. In recent decades, the long-time behavior of solutions to such equations has been extensively explored in the literature, see, e.g., \cite{jiao5,jiao3,jiao1,jiao,jiao2,jiao4} and the references therein. Moreover, the existence of attractors for  nonlocal  fractional deterministic differential equations has  been examined  in \cite{jiao6,jiao7}, whereas for the stochastic case, the existence of attractors has been investigated in \cite{jiao8,jiao9,jiao10,jiao11}. Furthermore, the existence of invariant measures for  nonlocal  fractional differential equations has been discussed in \cite{jiao12}. Meanwhile, stochastic evolution equations with memory have been studied extensively in the literature from different viewpoints.
For example, the pullback attractors to integro-differential equation with exponential-type memory have been established in \cite{caraballo2010}.
The  global attractors to stochastic evolution equations with memory have been investigated,
see, e.g., \cite{2,3,jiao13} and \cite{jiao14,jiao15} for deterministic and stochastic problems, respectively.

However, we note that invariant measures for stochastic evolution equations with memory have received much less attention.
 The main reason is that the process $u(t)$ governed by \eqref{1.1} is
not Markovian due to the presence of memory effects introduced by the convolution term
involving the kernel $K$. Consequently, it is not meaningful at this stage to speak of an
invariant probability measure for the process $u(t)$ alone. To overcome this difficulty,  we augment the original system by introducing a suitable
``history variable'' that records the past states of the solution. This idea goes back to \cite{jiao16} and has subsequently  been adopted  in \cite{3} for the study of singular limits in systems with memory. Later, the authors in \cite{er1} investigated the existence and regularity of invariant measures for the stochastic Coleman-Gurtin equation with additive noise. To be precise, we introduce the integrated past history of $u(t)$ defined by
\begin{align}\label{li1}
    \eta(t,s):= \int_{0}^{s} u(t-r)\mathrm{d}r=\int_{t-s}^{t} u(r)\mathrm{d}r,
    \qquad  s,t\geq0.
\end{align}
One easily checks that $\eta$ satisfies the inhomogeneous transport type equation,
\begin{align*}
    \partial_{t}\eta(t,s)=-\partial_{s}\eta(t,s)+u(t).
\end{align*}
To reveal the role of $\eta$ in  memory term in \eqref{1.1}, we define
\begin{align}\label{new10}
K(s)=\int_{s}^{\infty}\mu(r)\mathrm{d}r,\qquad s>0.
\end{align}
It means $\mu(s):=-K'(s)$ and $K(\infty)=0$.
Integrating by parts with respect to $s$ and using \eqref{li1}, we obtain
\[
\int_{0}^{\infty} K(s)(-\Delta)^{\alpha}u(t-s)\mathrm{d}s
    =-\int_{0}^{\infty}K'(s)(-\Delta)^{\alpha}\eta(t,s)\mathrm{d}s.
\]

This shows that the original non-Markovian equation \eqref{1.1} can be
rewritten as an extended system coupling the variables $(u(t),\eta(t,s))$, whose
dynamics are Markovian in an enlarged phase space.
The reformulation is essential in establishing well-posedness and  existence of
invariant measures for the associated Markov semigroup.
Thus, problem (1.1) is equivalent to the following extended system,
\begin{equation}\label{1.2}
\left\{
\begin{aligned}
&\mathrm{d}u(t)=
-\kappa(-\Delta)^{\alpha}u(t)\mathrm{d}t
-(1-\kappa)\int_{0}^{\infty}\mu(s)(-\Delta)^{\alpha}\eta(t,s)\mathrm{d}s\mathrm{d}t
+\varphi(u(t))\mathrm{d}t-\lambda u(t)\mathrm{d}t \\[0.8ex]
&+\sigma(u(t))\mathrm{d}W(t),\qquad (t,x)\in \mathbb{R}^{+}\times \mathcal{O},\\[0.8ex]
&\partial_t\eta(t,s)=-\partial_s\eta(t,s)+u(t),\qquad (t,s,x)\in \mathbb{R}^{+}\times\mathbb{R}^{+}\times \mathcal{O},\\[0.8ex]
&u(t)\big|_{\partial\mathcal{O}}=0,\qquad\eta(t,s)\big|_{\partial\mathcal{O}}=0,\qquad (t,s)\in \mathbb{R}^{+}\times \mathbb{R}^{+},\\[0.8ex]
&u(0)=u_{0},\quad\eta(0,s)=\eta_{0}(s):=\int_{-s}^{0}\psi(r)\mathrm{d}r,\quad (s,x)\in \mathbb{R}^+\times \mathcal{O}.
\end{aligned}
\right.
\end{equation}
We point out that, as far as the authors are aware, there are few  results concerning  invariant measures for stochastic partial differential equations  driven by nonlinear noise with memory. In the special case where the nonlinear noise reduces to the additive case, the authors in \cite{er1} investigated the regularity of invariant measures by means of a control argument.
 However, for stochastic fractional partial differential equations involving both nonlinear noise and memory terms, the well-posedness of solutions, together with the existence and regularity of invariant measures, remains an open problem.

Motivated by the above observations, we outline the three main contributions of this paper. {\bf (i)} The first objective is to establish the well-posedness of solutions to \eqref{1.2} with nonlinear multiplicative noise. In contrast to the additive noise case considered in \cite{er1}, the corresponding well-posedness result to the present nonlinear stochastic perturbation appears to be new. To this end, we reformulate the problem in an extended phase space and construct the solutions by means of Faedo-Galerkin approximations together with Minty's method. {\bf (ii)} The second objective is to prove the existence of invariant measures for the Markov semigroup generated by the system \eqref{1.2}. This is nontrivial because, although the introduction of the history variable restores the Markov property, the memory component prevents compactness in the natural phase space. To tackle this difficulty, we introduce suitable memory adapted functional spaces and combine them with uniform moment estimates to recover the compactness required by the Krylov-Bogoliubov argument. {\bf (iii)} The third objective is to investigate the regularity of invariant measures. To be precise, we show that the invariant measure is supported on a more regular phase space. The main observation is that the hyperbolic memory equation does not generate additional regularity, but only carries forward the regularity already available in the system. Higher regularity is therefore obtained through the smoothing effect of the fractional diffusion, together with a suitably chosen finite-dimensional control, which forces the long-time dynamics to enter a more regular phase space.

The paper is organized as follows. In Section 2, we present the functional framework, state the assumptions on the drift term, the diffusion coefficient and the memory kernel, and reformulate the original problem as a Markovian evolution equation in an extended phase space. Section 3 is devoted to the well-posedness of solutions, which is established by means of Faedo-Galerkin approximations and a priori estimates. In Section 4, we prove the existence of invariant measures for the corresponding Markov semigroup. In Section 5, we study the regularity of invariant measures by deriving higher-order estimates and applying a finite-dimensional control argument. Throughout this paper,  the letters $c$ and $C$  may denote different constants from line to line even in the same line.

\section{Preliminaries}
\setcounter{equation}{0}
In this section, we recall the functional settings and basic assumptions used throughout the paper and  gather several auxiliary results needed for
analysis of our problem.
\subsection{Fractional Laplacian}
Throughout this paper, we assume that $\mathcal O\subset\mathbb R^n$
is a bounded domain with smooth boundary, $0<\alpha<1$. Let $\mathcal{S}(\mathbb{R}^{n})$ denote the Schwartz space of smooth and rapidly decaying functions.
For $0<\alpha<1$, the fractional Laplacian $(-\Delta)^{\alpha}$ of a function
$u\in\mathcal{S}$ is defined by the singular integral,
\begin{align*}
(-\Delta)^{\alpha}u(x)
   =-\frac{1}{2}C(n,\alpha)
     \int_{\mathbb{R}^{n}}
     \frac{u(x+y)+u(x-y)-2u(x)}{|y|^{n+2\alpha}}\mathrm{d}y,
     \qquad x\in\mathbb{R}^{n},
\end{align*}
where the constant $C(n,\alpha)$ is positive and given by
\begin{align*}
    C(n,\alpha)
    = \frac{\alpha4^{\alpha}
      \Gamma\left(\tfrac{n+2\alpha}{2}\right)}
      {\pi^{\frac{n}{2}}\Gamma(1-\alpha)}.
\end{align*}
According to \cite{xin}, we find that for any $u\in\mathcal{S}$,
\begin{align*}
(-\Delta)^{\alpha}u
   =\mathcal{F}^{-1}\left(|\xi|^{2\alpha}\mathcal{F}u\right),
   \qquad\xi\in\mathbb{R}^{n},
\end{align*}
where $\mathcal{F}$ and $\mathcal{F}^{-1}$ denote the Fourier and inverse Fourier
transforms, respectively.

Let $H^{\alpha}(\mathbb{R}^{n})$ be the fractional Sobolev space defined by
\[
H^{\alpha}(\mathbb{R}^{n})
  = \left\{u\in L^{2}(\mathbb{R}^{n}):
      \iint_{\mathbb{R}^{n}\times\mathbb{R}^{n}}
      \frac{|u(x)-u(y)|^{2}}{|x-y|^{n+2\alpha}}\mathrm{d}x\mathrm{d}y <\infty
\right\},
\]
which is equipped with the norm
\[
\|u\|_{H^{\alpha}(\mathbb{R}^{n})}
   = \left(
       \int_{\mathbb{R}^{n}} |u(x)|^{2}\mathrm{d}x +\iint_{\mathbb{R}^{n}\times\mathbb{R}^{n}}
       \frac{|u(x)-u(y)|^{2}}{|x-y|^{n+2\alpha}}\mathrm{d}x\mathrm{d}y
     \right)^{1/2}.
\]
We also write the Gagliardo semi-norm of $H^{\alpha}(\mathbb{R}^{n})$ as
$\|\cdot\|_{\dot H^{\alpha}(\mathbb{R}^{n})}$, i.e.,
\begin{align}\label{liang}
\|u\|^{2}_{\dot H^{\alpha}(\mathbb{R}^{n})}
   = \int\int_{\mathbb{R}^{n}\times\mathbb{R}^{n}}
     \frac{|u(x)-u(y)|^{2}}{|x-y|^{\,n+2\alpha}}\mathrm{d}x\mathrm{d}y,
     \qquad u\in H^{\alpha}(\mathbb{R}^{n}).
\end{align}
Then, for all $u\in H^{\alpha}(\mathbb{R}^{n})$, we have
\[
\|u\|^{2}_{H^{\alpha}(\mathbb{R}^{n})}
   =\|u\|^{2}_{L^{2}(\mathbb{R}^{n})}+\|u\|^{2}_{\dot H^{\alpha}(\mathbb{R}^{n})}.
\]
Note that $H^{\alpha}(\mathbb{R}^{n})$ is a Hilbert space with inner product
\[
(u,v)_{H^{\alpha}(\mathbb{R}^{n})}
=\int_{\mathbb{R}^{n}} u(x)v(x)\mathrm{d}x
+\iint_{\mathbb{R}^{n}\times\mathbb{R}^{n}}
\frac{(u(x)-u(y))(v(x)-v(y))}{|x-y|^{n+2\alpha}}\mathrm{d}x\mathrm{d}y,
\qquad u,v\in H^{\alpha}(\mathbb{R}^{n}).
\]
It follows from  \cite{xin} that the norm of $\|u\|_{H^{\alpha}(\mathbb{R}^{n})}$ is equivalent to
\[
\|u\|^{2}_{H^{\alpha}(\mathbb{R}^{n})}
=\|u\|^{2}_{L^{2}(\mathbb{R}^{n})}
+\frac{2}{C(n,\alpha)}\|(-\Delta)^{\alpha/2}u\|^{2}_{L^{2}(\mathbb{R}^{n})},
\qquad u\in H^{\alpha}(\mathbb{R}^{n}),
\]
and hence
$\Big(\|u\|^{2}_{L^{2}(\mathbb{R}^{n})}+\|(-\Delta)^{\alpha/2}u\|^{2}_{L^{2}(\mathbb{R}^{n})}
\Big)^{1/2}$ is an equivalent norm of $H^{\alpha}(\mathbb{R}^{n})$.

Since the fractional Laplacian $(-\Delta)^{\alpha}$ is a
nonlocal operator, the homogeneous Dirichlet boundary condition \eqref{1.2}
is naturally interpreted as $u =0$ on $\mathbb{R}^{n}\setminus\mathcal{O}$
instead of  $u=0$ only on $\partial\mathcal{O}$.  Based on this interpretation, we set
$V=\{u\in H^{\alpha}(\mathbb{R}^{n}):u=0~\text{a.e.\ on }
\mathbb{R}^{n}\setminus\mathcal{O}\}$ and $H=\{u\in L^{2}(\mathbb{R}^{n}):u=0~\text{a.e.\ on }
\mathbb{R}^{n}\setminus\mathcal{O}\}$,
and denote by $\langle \cdot,\cdot\rangle_V$ and $\langle \cdot,\cdot\rangle_H$, $\|\cdot\|_V$ and $\|\cdot\|_H$
the inner products and norms in $V$ and $H$, respectively.
For $p>1$, we denote by $L^p(\mathcal{O})$ the usual $L^p$-space over $\mathcal{O}$
with the standard norm $\|\cdot\|_{L^{p}}$. Thus, we have the embeddings $V\hookrightarrow H=H^{*}\hookrightarrow V^{*},$
where $H^{*}$ and $V^{*}$ denote the dual spaces of $H$ and $V$, respectively,
and $H^{*}$ is identified with $H$ by Riesz's representation theorem. At last, we  use $\mathcal{L}_{2}(\mathcal{U},H)$ to denote the space of Hilbert--Schmidt
operators from a separable Hilbert space $\mathcal{U}$ to $H$ endowed with the norm $\|\cdot\|_{\mathcal{L}_{2}(\mathcal{U},H)}$.
Throughout this paper, we assume that
$(\Omega,\mathcal{F},\{\mathcal{F}_{t}\}_{t\in\mathbb R},\mathbb{P})$
is a complete filtered probability space satisfying the usual conditions.

Following the framework of Wang \cite{jiao10}, let $a_{\alpha}:V\times V\to\mathbb{R}$ be the bilinear form given by,
for $u,v\in V$,
\[a_{\alpha}(u,v)=\frac{C(n,\alpha)}{2}\int_{\mathbb{R}^{n}}\int_{\mathbb{R}^{n}}\frac{(u(x)-u(y))(v(x)-v(y))}{|x-y|^{n+2\alpha}}\mathrm{d}x\mathrm{d}y,\]
where $C(n,\alpha)$ is the constant appearing in the definition of the fractional Laplacian. For convenience, we associate an operator
$A_{\alpha}:V\to V^{*}$ with $a_{\alpha}$ such that
\[\langle A_{\alpha}u,v\rangle_{V^{*},V}=a_{\alpha}(u,v),~u,v\in V,\]
where $\langle\cdot,\cdot\rangle_{V^{*},V}$ denotes the duality pairing of $V^{*}$ and $V$.
 It is evident that the operator $A_{\alpha}$ is injective. On the other hand, by the Riesz representation theorem, $A_{\alpha}$ is also surjective and
thus the inverse $A_{\alpha}^{-1}:V^{*}\to V$ is well defined. Since $H=H^{*}$ can be considered as a subspace of
$V^{*}$ and the embedding $V\hookrightarrow H$ is compact, we find that $A_{\alpha}^{-1}:H\to V\subset H$ is a symmetric compact operator.
 It then follows from the Hilbert--Schmidt theorem that $A_{\alpha}$
has a family of eigenfunctions $\{l_j\}_{j=1}^{\infty}$ such that
$\{l_j\}_{j=1}^{\infty}$ forms an orthonormal basis of $H$. Moreover,
if $\lambda_j$ is the eigenvalue corresponding to $l_j$, i.e.,
\[A_{\alpha}l_j=\lambda_j l_j,\qquad j=1,2,\ldots,\]
then
\[0<\lambda_1\leq\lambda_2\leq\cdots\leq\lambda_j\to\infty\qquad\text{as }j\to\infty.\]
Moreover,  the operator $A_\alpha:V\to V^*$ introduced above can also be regarded as an operator on $H$,
with domain
\[D(A_\alpha)=\{u\in V:\ \exists f\in H,a_\alpha(u,v)=(f,v)_H,\forall v\in V\}.\]
For every $u\in D(A_\alpha)$, we define $A_\alpha u=f$.
In particular, $D(A_\alpha)\subset V\subset H.$
Given $u\in H$, we have
\[u=\sum_{j=1}^{\infty}a_j l_j,\qquad a_j=\int_{\mathbb{R}^n}u(x)l_j(x)\mathrm{d}x.\]
Then for $r>0$, define
\[A_{\alpha}^r u=\sum_{j=1}^{\infty}a_j\lambda_j^r l_j,\]
provided the series is convergent for $u\in H$. The domain of $A^r$ is denoted by $D(A_{\alpha}^r)$,
which is equipped with the norm$\|u\|_{D(A_{\alpha}^r)}=\|A_{\alpha}^r u\|_{H}$ for $u\in D(A^r)$.
For $k=1,\ldots,m,~k\in \mathbb{N}$, we set
\[\mathbb{H}^{k\alpha}:=D\!\left(A_\alpha^{k/2}\right),\langle u,v\rangle_{\mathbb{H}^{k\alpha}}
:=\langle A_\alpha^{k/2}u,A_\alpha^{k/2}v\rangle_H,~\|u\|_{\mathbb{H}^{k\alpha}}=\left\|A_\alpha^{k/2}u\right\|_H,\]
with the convention $\mathbb{H}^{0}:=H$.

\subsection{Setting up functional framework}
We will study the nonlocal fractional stochastic partial differential equation \eqref{1.2} when the boundary condition is replaced by $u = 0$ on $\mathbb{R}^{n}\setminus\mathcal{O}$. Namely, for $\alpha\in(0,1)$, we consider the following equation,
\begin{equation}\label{liran1}
\left\{
\begin{aligned}
&\mathrm{d}u(t)=
-\kappa(-\Delta)^{\alpha}u(t)\mathrm{d}t
-(1-\kappa)\int_{0}^{\infty}\mu(s)(-\Delta)^{\alpha}\eta(t,s)\mathrm{d}s\mathrm{d}t
+\varphi(u(t))\mathrm{d}t-\lambda u(t)\mathrm{d}t\\[0.8ex]
&+\sigma(u(t))\mathrm{d}W(t),\qquad (t,x)\in \mathbb{R}^{+}\times \mathcal{O},\\[0.8ex]
&\partial_t\eta(t,s)=-\partial_s\eta(t,s)+u(t),\qquad (t,s,x)\in \mathbb{R}^{+}\times \mathbb{R}^{+}\times \mathcal{O},\\[0.8ex]
&u(t,x)=0,\qquad\eta(t,s,x)=0,~~~ x\in\mathbb{R}^{n}\setminus\mathcal{O},\\[0.8ex]
&u(0)=u_{0},\qquad\eta(0,s)=\eta_{0}(s),\quad (s,x)\in\mathbb{R}^+\times \mathcal{O}.
\end{aligned}
\right.
\end{equation}
 Throughout this paper, we assume the nonlinearity $\varphi : \mathbb{R} \to \mathbb{R}$ is a continuous function such that for any $s\in\mathbb{R}$:
 \begin{enumerate}
    \item[($P_0$)] $\varphi\in C^{1}(\mathbb{R})$ satisfying $\varphi(0)=0$;
    \item[($P_1$)] There exist positive constants $a_{1}>0$, $p_{0}\geq 1$, such that
    \[ |\varphi(s)|\leq a_{1}\bigl(1+|s|^{p_{0}}\bigr);\]
    \item[$(P_2)$] There exist positive constants $a_{2}, a_{3}$, such that
    \[s\varphi(s)\leq-a_{2}|s|^{p_{0}+1}+a_{3};\]
    \item[($P_3)$] There exists a constant $a_{\varphi}$, such that
    \[\sup_{x\in\mathbb{R}}\varphi'(x)=: a_{\varphi}<\infty.\]
\end{enumerate}
%where $p_{0}+1 \in \left(2,\frac{2n}{n-2\alpha}\right)$ such that
%$V$ is continuously embedded in
%$L^{p_{0}+1}(\mathcal{O})$.

Additionally, assume the diffusion coefficient $\sigma:H \to \mathcal{L}_{2}(\mathcal{U},H)$ is a nonlinear Lipschitz mapping fulfilling:
\begin{enumerate}
\item[($Q_0$)]
There exists a constant $L_{\sigma}>0$ such that, for all $u,v\in H$,
\[\|\sigma(u)-\sigma(v)\|_{\mathcal{L}_{2}(\mathcal{U},H)}\leq L_\sigma\|u-v\|_{H};\]
\end{enumerate}
From ($Q_0$), the following condition holds:
\begin{enumerate}
\item[($Q_1$)]
There exists a constant $C_{\sigma}>0$ such that for all $u\in H$,
\[\|\sigma(u)\|_{\mathcal{L}_{2}(\mathcal{U},H)}\leq C_{\sigma}(1+\|u\|_{H}),\]
where $C_{\sigma}=\max\{ \|\sigma(0)\|_{\mathcal{L}_{2}(\mathcal{U},H)}, L_\sigma\}.$
\end{enumerate}
\subsection{Time-dependent memory spaces}
In what follows, we introduce the concept of extended phase spaces that were developed in
\cite{2,3}. Let the memory kernel $\mu\in C^{1}\big([0,\infty)\big)$ be a positive function, such that
\begin{equation}\label{lii2}
\mu'+\delta\mu\leq0,
\quad\text{for some }\delta>0.
\end{equation}
We define the weighted Hilbert spaces
\[M_{\mu}^{0}:=L_{\mu}^{2}([0,\infty);V),\quad M_{\mu}^{k\alpha}:=L_{\mu}^{2}\left([0,\infty);D(A_{\alpha}^{(k+1)/2})\right),\quad k=1,2,\ldots,m,\]
equipped, respectively, with the inner products
\[\langle \eta_1,\eta_2\rangle_{M_{\mu}^{0}}=\int_0^\infty \mu(s)a_{\alpha}(\eta_1(s),\eta_2(s))\mathrm{d}s,\]
and
$$\langle \eta_1,\eta_2\rangle_{M_{\mu}^{k\alpha}}=\int_0^\infty \mu(s)\langle A_{\alpha}^{(k+1)/2}\eta_1(s),
A_{\alpha}^{(k+1)/2}\eta_2(s)\rangle_H\mathrm{d}s,\quad k=1,2,\ldots,m.$$
Following the approach of Xu et al. \cite{you}, we introduce the weighted past-history space
$$L_\gamma^2((-\infty,0);V)=\left\{\psi:(-\infty,0)\to V:\int_{-\infty}^0 e^{\gamma r}\|\psi(r)\|_V^2\mathrm{d}r<\infty\right\},$$
where $0<\gamma<\delta$.
Define the operator $J:L_\gamma^2((-\infty,0);V)\longrightarrow M_\mu^0$, satisfying
$$J\psi(s)=\int_{-s}^0\psi(r)\mathrm{d}r,\qquad s\geq0.$$
\begin{lem}\label{you}
Assume that the memory kernel $\mu$ satisfies \eqref{lii2}. Then the operator
$$J:L_\gamma^2((-\infty,0);V)\longrightarrow M_\mu^0$$
is bounded linear operator.
\end{lem}
\begin{proof}
The linearity of $J$ follows from the linearity of integration. Thus, it remains to establish the boundedness of $J$. Let $\psi\in L_\gamma^2((-\infty,0);V)$ be arbitrary. For any $s>0$, applying the Cauchy-Schwarz inequality, we have
\begin{align*}
\|J\psi(s)\|_V^2=\left\|\int_{-s}^0\psi(r)\mathrm{d}r\right\|_V^2\leq
\left(\int_{-s}^0e^{-\gamma r}\mathrm{d}r\right)\left(\int_{-s}^0e^{\gamma r}\|\psi(r)\|_V^2\mathrm{d}r\right)\leq\frac{e^{\gamma s}-1}{\gamma}\|\psi\|_{L_\gamma^2((-\infty,0);V)}^2.
\end{align*}
Consequently,
$$\|J\psi\|_{M_\mu^0}^2=\int_0^\infty\mu(s)\|J\psi(s)\|_V^2\mathrm{d}s\leq\frac1\gamma\left(\int_0^\infty\mu(s)(e^{\gamma s}-1)\mathrm{d}s\right)\|\psi\|_{L_\gamma^2((-\infty,0);V)}^2.$$
Since \eqref{lii2} implies $\mu(s)\leq\mu(0)e^{-\delta s}$ and $0<\gamma<\delta$, we know
$$\|J\psi\|_{M_\mu^0}^2\leq\frac{\mu(0)}{\delta(\delta-\gamma)}\|\psi\|_{L_\gamma^2((-\infty,0);V)}^2.$$
Hence $J$ is a bounded linear operator.
\end{proof}
In particular, for every $\psi\in L_\gamma^2((-\infty,0);V)$, the initial memory variable $\eta_0(s)=\int_{-s}^0\psi(r)\mathrm{d}r$ belongs to $ M_\mu^0$.

The embedding
$M_{\mu}^{k_1\alpha}\hookrightarrow M_{\mu}^{k_2\alpha}$ is merely continuous whenever $k_1>k_2$, with $k_1,k_2=0,1,\ldots,m$.
Therefore, it is necessary to introduce additional functional spaces constructed below to establish the existence of invariant measures  by the classical Krylov-Bogoliubov approach.

Let $\mathcal{T}_{\mu}$ be the operator on $M^{0}_{\mu}$ defined by
\begin{align*}
\mathcal{T}_{\mu}\eta:=-\partial_{s}\eta,
\qquad \mathrm{Dom}(\mathcal{T}_{\mu}) =\{\eta\in M^{0}_{\mu}:\partial_{s}\eta\in M^{0}_{\mu},~\eta(0)=0\},
\end{align*}
where $\partial_{s}$ is the derivative with respect to $s$ in the distribution sense.
In fact, $\mathcal{T}_{\mu}$ is the infinitesimal generator of the right-translation semigroup
acting on $M^{0}_{\mu}$~\cite[Theorem~3.1]{6}. Furthermore, we can obtain the following lemma.
\begin{lem}\label{mon}
If $u\in L_{\mathrm{loc}}^{1}([0,\infty);V)$ , then the following functional formulation
of the Cauchy initial value problem
\begin{align}\label{nan}
\begin{cases}
\dfrac{\mathrm{d}}{\mathrm{d}t}\eta(t)=\mathcal{T}_{\mu}\eta(t)+u(t),\\[1.5ex]
\eta(0,s)=\eta_{0}(s)\in M_{\mu}^{0},
\end{cases}
\end{align}
has a unique solution and this solution, corresponding to the initial data $\eta_{0}$, can be expressed as
$\eta\in C\big([0,\infty);M^{0}_{\mu}\big),$
\begin{align}\label{x}
\eta(t,s;\eta_{0}) =
\begin{cases}
\displaystyle \int_{0}^{s}u(t-r)\mathrm{d}r,
&0<s\leq t, \\
\displaystyle \eta_{0}(s-t)+\int_{0}^{t}u(t-r)\mathrm{d}r,
&s>t.
\end{cases}
\end{align}
\end{lem}

Given $\eta \in M_{\mu}^{0}$, we next introduce the \emph{tail function}
\begin{equation}\label{2.9}
\mathbb{T}_{\eta}^{\mu}(r)=\int_{(0,\frac{1}{r})\cup(r,\infty)}\mu(s)
a_{\alpha}(\eta(s),\eta(s))\mathrm{d}s,\qquad r\geq1.
\end{equation}
For $k=0,1,\ldots,m$, we define the Banach space
\begin{align*}\label{lii5}
E_{\mu}^{k\alpha}
=\left\{\eta\in M_{\mu}^{k\alpha}:\eta\in\mathrm{Dom}(\mathcal{T}_{\mu}), ~\sup_{r\geq1}r \mathbb{T}_{\eta}^{\mu}(r)<\infty\right\}.
\end{align*}
The associated norm is given by
\begin{equation}\label{2.11}
\|\eta\|_{E_{\mu}^{k\alpha}}^{2}=\|\eta\|_{M_{\mu}^{k\alpha}}^{2}+\|\mathcal{T}_{\mu}\eta\|_{M_{\mu}^{0}}^{2}
+\sup_{r\geq1}r\mathbb{T}_{\eta}^{\mu}(r).
\end{equation}
Having established the memory spaces, we next define the Banach spaces associated with any
$k=0,1,\ldots,m$,
\[\mathcal H_\mu^{k\alpha}=\mathbb{H}^{k\alpha}\times M_\mu^{k\alpha},~\mathcal{Z}_\mu^{k\alpha}=\mathbb{H}^{k\alpha}\times E_\mu^{k\alpha},~k=0,1,\ldots,m.\]
For $k=1,\ldots,m$, $E_{\mu}^{k\alpha}$ is compactly embedded into $M_{\mu}^{0}$ \cite[Lemma~5.5]{5},
so that any bounded set in $\mathcal{Z}_{\mu}^{k\alpha}$ is precompact in $\mathcal{H}_{\mu}^{0}$.
%\green{Furthermore, for $n\geq1$, we denote,
%\begin{align}\label{chen2}
%P_{n}u=\sum_{k=1}^{n}\langle u,e_{k}\rangle_{H}e_{k},
%\qquad P_{n}\eta(s)=\sum_{k=1}^{n}\langle\eta(s),e_{k}\rangle_{H}e_{k},
%\end{align}
%where $\{e_k\}_{k\geq1}$ denotes an orthonormal basis of the Hilbert space $H$. The operator $P_n$ is the orthogonal projection from $H$ onto the finite-dimensional subspace
%$H_n := \mathrm{span}\{e_1,\ldots,e_n\}$.}

In addition, the operator \(\mathcal{T}_{\mu}\) has the
following estimate (see \cite[Theorem~3.1]{6}), which will be
frequently used in this  manuscript.
For any \(\eta \in \mathrm{Dom}(\mathcal{T}_{\mu})\), using integration by parts and \eqref{lii2}, we obtain
\begin{equation}\label{lii1}
\begin{aligned}
\langle\mathcal{T}_{\mu}\eta,\eta\rangle_{M^{0}_{\mu}}
&=-\frac12\int_{0}^{\infty}\mu(s)\partial_{s}a_{\alpha}(\eta(s),\eta(s))\mathrm{d}s\\[0.8ex]
&=\frac12\int_{0}^{\infty}\mu'(s)a_{\alpha}(\eta(s),\eta(s))\mathrm{d}s\\[0.8ex]
&\leq-\frac12\delta\|\eta\|^{2}_{M^{0}_{\mu}}.
\end{aligned}
\end{equation}
In a similar manner, we derive the same estimate that for any $k=1,\ldots,m$,
\begin{equation}\label{2.8}
\langle\mathcal{T}_{\mu}\eta,\eta\rangle_{M_{\mu}^{k\alpha}}
\leq-\frac{1}{2}\delta\|\eta\|^{2}_{M_{\mu}^{k\alpha}}.
\end{equation}
The above formula will play a key role in the proof of the
existence of invariant measures.

Based on the above construction, we can now rewrite problem \eqref{liran1} in the following form:
\begin{equation}\label{lii10}
\left\{
\begin{aligned}
&\mathrm{d}u(t)=
-\kappa A_{\alpha}u(t)\mathrm{d}t
-(1-\kappa)\int_{0}^{\infty}\mu(s)A_{\alpha}\eta(t,s)\mathrm{d}s\mathrm{d}t
+\varphi(u(t))\mathrm{d}t-\lambda u(t)\mathrm{d}t\\[0.8ex]
&+\sigma(u(t))\mathrm{d}W(t),\qquad (t,x)\in \mathbb{R}^{+}\times \mathcal{O},\\[0.8ex]
&\frac{\mathrm{d}}{\mathrm{d}t}\eta(t)=\mathcal{T}_\mu\eta(t)+u(t),\quad (t,x)\in\mathbb{R}^+\times \mathcal{O},\\[0.8ex]
&\left(u(t,x),\eta(t,s,x)\right)=0,~~~ x\in\mathbb{R}^{n}\setminus\mathcal{O},\\[0.8ex]
&u(0)=u_{0},\qquad\eta(0,s)=\eta_{0}(s),\quad (s,x)\in\mathbb{R}^+\times \mathcal{O}.
\end{aligned}
\right.
\end{equation}

We denote by $\mathcal{L}^\mu$ the generator corresponding to system \eqref{lii10}.
For any function $g \in C^{2}(\mathcal{H}^0_\mu)$ of the form $g=g(U)$ with $U=(u,\eta)$,
the operator $\mathcal{L}^\mu$ is defined by
\begin{equation}\label{lili}
\begin{aligned}
\mathcal{L}^\mu g(u,\eta):=&-\kappa\langle A_{\alpha}u,D_u g\rangle_H
-(1-\kappa)\int_{0}^{\infty}\mu(s)\langle A_{\alpha}\eta(s), D_u g\rangle_H \mathrm{d}s-\langle\lambda u,D_u g\rangle_H
\\
&+\langle \varphi(u), D_u g\rangle_{\langle L^{q}(\mathcal{O}),L^{p_{0}+1}(\mathcal{O})\rangle}+\langle \mathcal{T}_\mu\eta, D_\eta g\rangle_{M_\mu^0}
+\langle u,D_\eta g\rangle_{M_\mu^0}
+\frac12\operatorname{Tr}\left(D_{uu} g\, \sigma\sigma^{*}\right),
\end{aligned}
\end{equation}
where all the terms above are well defined.
\section{Well-posedness}
\setcounter{equation}{0}
This section is devoted to establishing  the existence and uniqueness of solutions to problem \eqref{lii10} in the following sense.
\begin{definition}\label{dingyi1}
Fix $T>0$ and let
$U_{0}=(u_0,\eta_0)\in \mathcal H_{\mu}^{0}$. Then, a continuous   process
$U(\cdot)=(u(\cdot),\eta(\cdot))$ is called a strong solution (in the probabilistic sense) of \eqref{lii10} on $[0,T]$, if $U(\cdot)$ is $\mathcal{F}_t$-adapted fulfilling
\begin{align*}
&u\in L^{2}(\Omega;C([0,T];H))\cap L^{2}(\Omega;L^2([0,T];V))
\cap L^{p_0+1}(\Omega;L^{p_0+1}([0,T];L^{p_0+1}(\mathcal{O}))),\\
&\eta\in L^{2}(\Omega;C([0,T];M^0_\mu)),
\end{align*}
such that, for all $v \in V \cap L^{p_{0}+1}(\mathcal{O})$ and $\zeta\in M^0_\mu$,
\begin{align*}
&\langle u(t),v\rangle_H+\kappa\int_0^t a_{\alpha}( u(r),
v)\mathrm{d}r+\lambda\int_0^t\langle u(r),v\rangle_{H}\mathrm{d}r+(1-\kappa)\int_0^t\int_0^\infty\mu(l)
a_{\alpha}(\eta(r,l),v)\mathrm{d}l\mathrm{d}r\notag\\[0.8ex]
&-\int_0^t\langle\varphi(u(r)),v\rangle_{\langle L^{q}(\mathcal{O}),L^{p_{0}+1}(\mathcal{O})\rangle}\mathrm{d}r\notag\\[0.8ex]
=&\langle u_0,v\rangle_H+\int_0^t\langle\sigma(u(r))\mathrm{d}W(r),v\rangle_H,
\end{align*}
where  $q$  is the H\"older conjugate of $p_{0}+1$ and
\begin{align*}
\langle \eta(t),\zeta\rangle_{M^0_\mu}=\langle\eta_0,\zeta\rangle_{M^0_\mu}
+\int_0^t\langle \mathcal{T}_\mu\eta(r),\zeta\rangle_{M^0_\mu}\mathrm{d}r
+\int_0^t\langle u(r),\zeta\rangle_{M^0_\mu}\mathrm{d}r,
\end{align*}
hold $\mathbb P$-almost surely.
\end{definition}
\setcounter{theorem}{1}
\begin{theorem}\label{the1}
Assume that $(Q_0)$ and $(P_0)$-$(P_3)$ hold.
Then, for all $U_0:=(u_{0},\eta_{0})\in\mathcal{H}^0_{\mu}$ and  $T>0$, problem \eqref{lii10}  admits
a unique strong solution $U(\cdot;U_0)$ in the sense of Definition \ref{dingyi1}. Moreover, the solution depends continuously on the initial data $U_0$ in $\mathcal H_{\mu}^{0}$.
\end{theorem}
\begin{proof}
We prove the well-posedness of \eqref{lii10} by the classical Faedo-Galerkin approximations and Minty's method for monotone operators \cite{NEW}. The proof consists of the following steps.

\textbf{Step 1. Approximate solutions.}
We construct a sequence of approximate solutions using finite-dimensional Galerkin approximations. Since $V\cap L^{p_0+1}(\mathcal O)$ is a separable
Banach space densely and continuously embedded into $H$, there exists an orthonormal basis $\{e_j\}_{j=1}^{\infty}$ of $H$ such that
$\{e_j\}_{j=1}^{\infty}\subset V\cap L^{p_0+1}(\mathcal O)$ and $\operatorname{span}\{e_j:j\in\mathbb N\}$ is dense in $V\cap L^{p_0+1}(\mathcal O)$.
For $N\in\mathbb N$, let $H_N:=\operatorname{span}\{e_1,\ldots,e_N\}$, and define the orthogonal projection $P_N:H\to H_N$ by
\begin{equation}
P_Nu=\sum_{j=1}^{N}\langle u,e_j\rangle_H e_j,
\qquad u\in H.
\end{equation}
For the memory component, we have
\[(P_N\eta)(s):=\sum_{j=1}^{N}\langle\eta(s),e_j\rangle_H e_j,\qquad \eta\in M_\mu^0,\]
for almost every $s\geq 0$.
In fact, we can extend $P_N$ to $V^*$ and
$L^q(\mathcal{O})$ by
\[P_N\phi=\sum_{j=1}^{n}\phi(e_j)e_j,\qquad\text{for } \phi\in V^*\text{ or } \phi\in L^q(\mathcal{O}).\]
Let $\{g_j\}_{j=1}^{\infty}$ be an orthonormal basis of $\mathcal{U}$, and let $\{\beta_j\}_{j=1}^{\infty}$ be the corresponding independent standard Brownian motions associated with $W$.
For $N\in\mathbb N$, set
\[\mathcal{U}_N:=\operatorname{span}\{g_1,\ldots,g_N\},\qquad W_N(t):=\sum_{j=1}^N\beta_j(t)g_j.\]
We seek approximate solutions to problem \eqref{lii10} of the form
\[u_N(t,x)=\sum_{j=1}^N u_j^N(t)e_j(x),\qquad\eta_N(t,s,x)=\sum_{j=1}^N\eta_j^N(t,s)e_j(x).\]
Then, we introduce the orthogonal projection $R_N:V\longrightarrow H_N$
with respect to the bilinear form $a_\alpha$, satisfying
\[a_\alpha(R_Nv,w)=a_\alpha(v,w),\qquad \forall v\in V,\quad \forall w\in H_N.\]
Substituting the above approximations into the original system and taking the inner product with each basis
function $e_j$, we obtain the following finite-dimensional system,
\begin{equation}\label{lili1}
\begin{cases}
\displaystyle
\mathrm{d}u_N(t)=-\kappa P_{N}A_{\alpha}u_N(t)\mathrm{d}t-\lambda u_N(t)\mathrm{d}t
-(1-\kappa)P_{N}\int_0^\infty \mu(s)A_{\alpha}\eta_N(t,s)\mathrm{d}s\mathrm{d}t
+P_N\varphi(u_N(t))\mathrm{d}t\\[0.8ex]
+P_N\sigma(u_N(t))\mathrm{d}W_{N}(t),\qquad (t,x)\in \mathbb{R}^{+}\times \mathcal{O},\\[0.8ex]
\frac{\mathrm{d}}{\mathrm{d}t}\eta_N(t)=\mathcal{T}_\mu\eta_N(t)+u_N(t),\qquad (t,x)\in \mathbb{R}^{+}\times \mathcal{O}\\[0.8ex]
u_N(0)=P_Nu_{0},\qquad \eta_N(0,s)=R_N\eta_{0}(s),\qquad\eta_N(t,0)=0.
\end{cases}
\end{equation}
According to standard existence theory for SODE there exists a continuous solution
of \eqref{lili1} on some interval $(0,T)$. The a priori estimates that follow
imply that, in fact, $T=+\infty$.

\textbf{Step 2. Uniform estimates.}
We will derive uniform estimates of approximate solutions to problem \eqref{lili1} with respect to $N \in \mathbb{N}$. To this end, let us define  $ U_N(t)=(u_N(t),\eta_N(t))$ and introduce the function $\Psi_0(u,\eta)$ given by
\begin{align}\label{shu}
\Psi_0(u,\eta)=\frac12 \|u\|_H^2+\frac{1-\kappa}{2}\|\eta\|^2_{M_\mu^0}.
\end{align}
Computing partial derivatives of $\Psi_0$, one has
\[D_u\Psi_0=u,\quad D_\eta\Psi_0=(1-\kappa)\eta,\quad\text{and}\quad D_{uu}\Psi_0=Id_H.\]
Using It\^{o}'s formula to $\Psi_{0},$ replacing $(u,\eta)$ by $(u_{N},\eta_{N})$ and combining with \eqref{lili1}, we obtain
\begin{align}\label{1.4}
\mathrm{d}\Psi_0(U_N(t))
=&-\frac{\kappa C(n,\alpha)}{2}\|u_N(t)\| _{\dot{H}^{\alpha}(\mathbb{R}^{n})}^2\mathrm{d}t-(1-\kappa)\langle u_N(t),\eta_N(t)\rangle_{M_\mu^{0}}\mathrm{d}t-\lambda\|u_N(t)\|_{H}^{2}\mathrm{d}t\notag\\
&+\langle P_N\varphi(u_N(t)),u_N(t)\rangle_{\langle L^{q}(\mathcal{O}),L^{p_{0}+1}(\mathcal{O})\rangle}\mathrm{d}t
+(1-\kappa)\langle\eta_N(t),\mathcal{T}_\mu\eta_N(t)\rangle_{M_\mu^0}\mathrm{d}t\notag\\
&+(1-\kappa)\langle\eta_N(t), u_N(t)\rangle_{M_\mu^0}\mathrm{d}t
+\frac12\|P_N \sigma(u_N(t))\|_{\mathcal{L}_2(\mathcal{U}_{N},H)}^2\mathrm{d}t\notag\\
&+\langle u_N(t), P_N\sigma(u_N(t))\mathrm{d}W_{N}(t)\rangle_H.
\end{align}
Since $u_N\in H_N$ and $P_N$ is extended to $L^q(\mathcal O)$ by duality, by assumption $(P_2)$, we have
\[\langle P_N\varphi(u_N),u_N\rangle_{\langle L^{q}(\mathcal{O}),L^{p_{0}+1}(\mathcal{O})\rangle}=\langle\varphi(u_N),u_N\rangle_{\langle L^{q}(\mathcal{O}),L^{p_{0}+1}(\mathcal{O})\rangle}
\leq-a_2\|u_N\|_{L^{p_0+1}}^{p_0+1}+a_3|\mathcal{O}|.\]
By the linear growth of \(\sigma\) in $(Q_1)$, we infer
\[\|P_N \sigma(u_N)\|_{\mathcal{L}_2(\mathcal{U}_{N},H)}^2\leq \|\sigma(u_N)\|_{\mathcal{L}_2(\mathcal{U}_{N},H)}^2\leq \|\sigma(u_N)\|_{\mathcal{L}_2(\mathcal{U},H)}^2
\leq 2C_{\sigma}^{2}\big(1+\|u_N\|_H^2\big).\]
It follows from the above estimates and \eqref{lii1} that
\begin{align}\label{1.5}
&\mathrm{d}\Psi_0(U_N(t))+\frac{\kappa C(n,\alpha)}{2}\|u_N(t)\| _{\dot{H}^{\alpha}(\mathbb{R}^{n})}^2\mathrm{d}t
+\frac{1}{2}(1-\kappa)\delta\|\eta_N(t)\|_{M_\mu^0}^2\mathrm{d}t\notag\\[0.6ex]
& +\lambda\|u_N(t)\|_H^2\mathrm{d}t+a_2\|u_N(t)\|_{L^{p_0+1}}^{p_0+1}\mathrm{d}t\notag\\[1.0ex]
&\leq C\big(1+\|u_N(t)\|_H^2\big)\mathrm{d}t+\langle u_N(t),P_N\sigma(u_N(t))\mathrm{d}W_{N}(t)\rangle_H,
\end{align}
where $C$ is dependent on $a_3$, $\mathcal{O}$ and $C_{\sigma}$. Integrating \eqref{1.5} from 0 to $t$ and then taking the supremum with respect to $t$ and expectation, we derive
\begin{align}\label{1.6}
&\mathbb{E}\left(\sup_{0\leq s\leq t}\Psi_0(U_N(s))\right)+\frac{\kappa C(n,\alpha)}{2}\mathbb{E}\int_0^t \|u_N(s)\|_{\dot{H}^{\alpha}(\mathbb{R}^{n})}^2\mathrm{d}s
+\frac12(1-\kappa)\delta\mathbb{E}\int_0^t\|\eta_N(s)\|_{M_\mu^0}^2\mathrm{d}s\notag\\
&+\lambda\int_0^t \mathbb{E}\|u_N(s)\|_H^2\mathrm{d}s+a_2\mathbb{E}\int_0^t\|u_N(s)\|_{L^{p_0+1}}^{p_0+1}\mathrm{d}s\notag\\
\leq&2\Psi_0(P_{N}u_{0},R_{N}\eta_{0})+2Ct+2C\int_0^t \mathbb{E}\left(\sup_{0\leq r\leq s}\|u_N(r)\|_H^2\right)\mathrm{d}s\notag\\
&+2\mathbb{E}\left(\sup_{0\leq s\leq t}\left|\int_0^s\langle u_N(r),P_N\sigma(u_N(r))\mathrm{d}W_{N}(r)\rangle_H\right|\right).
\end{align}
For the last term, applying the BDG inequality and condition $(Q_1)$, we obtain
\begin{align}\label{1.7}
&\mathbb{E}\left(\sup_{0\leq s\leq t}\left|\int_0^s\langle u_N(r),P_N\sigma(u_N(r))\mathrm{d}W_{N}(r)\rangle_H\right|\right)\notag\\
\leq& C_{b}\mathbb{E}\left(\int_0^t \|u_N(s)\|_H^2\|P_N\sigma(u_N(s))\|_{\mathcal{L}_2(\mathcal{U}_{N},H)}^2\mathrm{d}s\right)^{1/2}\notag\\
%\leq&\frac18\mathbb{E}\left(\sup_{0\leq s\leq t}\|u_N(s)\|_H^2\right)+C\mathbb{E}\int_0^t\bigl(1+\|u_N(s)\|_H^2\bigr)\mathrm{d}s,\notag\\
\leq&\frac18\mathbb{E}\left(\sup_{0\leq s\leq t}\|u_N(s)\|_H^2\right)+C\int_0^t\bigl(1+\mathbb{E}(\sup_{0\leq r\leq s}\|u_N(r)\|_H^2)\bigr)\mathrm{d}s\notag\\
\leq&\frac14\mathbb{E}\left(\sup_{0\leq s\leq t}\Psi_0(U_N(s))\right)+C\int_0^t\bigl(1+\mathbb{E}(\sup_{0\leq r\leq s}\|u_N(r)\|_H^2)\bigr)\mathrm{d}s,
\end{align}
where $C_{b}$ is a constant associated with the BDG inequality and $C$ depends on $C_{b}$ and $C_{\sigma}$. By \eqref{1.6}-\eqref{1.7} and the Gronwall inequality,  there exists a
constant $C>0$, such that, for all $0<T^\prime<T$,
\begin{align}\label{1.8}
&\mathbb{E}\left(\sup_{0\leq t\leq T^\prime}\Psi_0(U_N(t))\right)
+\frac{\kappa C(n,\alpha)}{2}\mathbb{E}\int_0^{T^{\prime}} \|u_N(s)\|_{\dot{H}^{\alpha}(\mathbb{R}^{n})}^2\mathrm{d}s+\lambda\int_0^{T^\prime} \mathbb{E}\|u_N(s)\|_H^2\mathrm{d}s\notag\\
&+\frac12(1-\kappa)\delta\mathbb{E}\int_0^{T^\prime}\|\eta_N(t)\|_{M_\mu^0}^2\mathrm{d}t
+a_{2}\mathbb{E}\int_0^{T^\prime}\|u_N(t)\|_{L^{p_0+1}}^{p_0+1}\mathrm{d}t\notag\\[1.0ex]
\leq&Ce^{CT}\Big(1+\|u_0\|_H^{2}+\|\eta_0\|_{M_\mu^0}^{2}\Big).
\end{align}

It  thereby follows from \eqref{1.8}  that $u_N$ is bounded in \[
L^{2}(\Omega; L^{\infty}(0,T; H))\cap L^{2}(\Omega; L^{2}(0,T;V))\cap L^{p_{0}+1}(\Omega;L^{p_{0}+1}(0,T;L^{p_{0}+1}(\mathcal{O}))),\]
and $\eta_{N}$ is bounded in
\[L^2(\Omega;L^\infty(0,T;M_\mu^0)).\]
Furthermore, thanks to $(P_1)$  and $q=(p_0+1)^\prime$, we obtain
\begin{align}\label{lili2}
\mathbb{E}\int_{0}^{t}\|\varphi(u_N(r))\|_{L^{q}}^{q}\mathrm{d}r\leq&
\mathbb{E}\int_{0}^{t}\int_{\mathcal{O}}
\left(a_1 (1+|u_N(r)|^{p_0}) \right)^{q}\mathrm{d}x\mathrm{d}r\notag\\
\leq& ct+c\mathbb{E}\int_{0}^{t}\|u_N(r)\|_{L^{p_0+1}}^{p_0+1}\mathrm{d}r.
\end{align}
\textbf{Step 3: Existence of solutions.}
We prove the existence of solutions of problem \eqref{lii10} by considering the limiting
process of solutions of problem \eqref{lili1} as $N\rightarrow\infty$.
From \eqref{1.8} and \eqref{lili2}, we deduce that, up to a subsequence (still denoted by $(u_{N},\eta_{N})$),
\begin{align}\label{zhuo0}
u_N&\rightharpoonup u\quad
\text{weakly in }L^2(\Omega;L^2(0,T;V));\\
u_N&\overset{*}{\rightharpoonup} u\quad
\text{weakly-star in }L^2(\Omega;L^\infty(0,T;H));\\
u_N&\rightharpoonup u\quad
\text{weakly in }L^{p_0+1}(\Omega;L^{p_0+1}(0,T;L^{p_0+1}(\mathcal{O})));\\
\eta_N&\overset{*}{\rightharpoonup}\eta\quad
\text{weakly-star in }L^2(\Omega;L^\infty(0,T;M^0_\mu));\\
\varphi(u_N)&\rightharpoonup\chi\quad\text{weakly in }L^q(\Omega;L^q(0,T;L^q(\mathcal{O}))).
%P_N\sigma(u_N)&\rightharpoonup\chi_{1}\quad
%\text{weakly in }L^2\left(\Omega;L^2(0,T;\mathcal{L}_{2}(\mathcal{U},H))\right);\\
%P_N\sigma(r,u_N(r))\mathrm{d}W(r)&\rightharpoonup\chi_{1}(r)\mathrm{d}W(r)\quad
%\text{weakly in }L^2\left(\Omega;L^2(0,T;H)\right).
\end{align}
%Furthermore, we have
%\begin{align}\label{jia1}
%P_N\varphi(u_N)\rightharpoonup\chi\text{ weakly in }L^q(\Omega;L^q(0,T;L^q(\mathcal{O}))).
%\end{align}
In addition, let $\bar{\sigma}_N(u_N)$ denote the zero extension of
$P_N\sigma(u_N):\mathcal{U}_N\to H_N$ to $\mathcal{U}$, namely,
\[\bar{\sigma}_N(u_N)h=
\begin{cases}
P_N\sigma(u_N)h, & h\in \mathcal{U}_N,\\
0, & h\in \mathcal{U}_N^\perp.
\end{cases}
\]
Then $\|\bar{\sigma}_N(u_N)\|_{L_2(\mathcal{U},H)}^2=\|P_N\sigma(u_N)\|_{L_2(\mathcal{U}_N,H)}^2.$
By $(Q_1)$ and the contractivity of $P_N$, the sequence
$\{\bar{\sigma}_N(u_N(\cdot))\}_{N\in\mathbb N}$ is bounded in $L^2\bigl(\Omega;L^2(0,T;\mathcal{L}_2(\mathcal{U},H))\bigr).$
Hence, there exist a subsequence, still denoted by $\{\bar{\sigma}_N(u_N(\cdot))\}$ and an element
$\chi_1\in L^2\bigl(\Omega;L^2(0,T;\mathcal{L}_2(\mathcal{U},H))\bigr)$ such that
\begin{align}\label{3.17}
\bar{\sigma}_N(u_N(\cdot))\rightharpoonup\chi_1\quad\text{weakly in }L^2\bigl(\Omega;L^2(0,T;\mathcal{L}_2(\mathcal{U},H))\bigr).
\end{align}
Moreover, by $(Q_1)$ and the uniform estimates of $u_N$, the sequence $\{\sigma(u_N)\}_{N\in\mathbb{N}}$ is bounded in $L^2\bigl(\Omega;L^2(0,T;\mathcal{L}_2(\mathcal{U},H))\bigr).$
Then, let $\Phi\in L^2\bigl(\Omega;L^2(0,T;L_2(\mathcal{U},H))\bigr)$ be arbitrary,  by the convergence from $P_N\to I$ strongly in $H$ and the Hilbert-Schmidt property of $\Phi$, we obtain
\[\begin{aligned}
&\left|\left\langle\sigma(u_N)-\bar\sigma_N(u_N),\Phi\right\rangle\right|\\
\leq& C\left[\mathbb E\int_0^T\left(\sum_{j=1}^N\|(I-P_N)\Phi g_j\|_H^2+\sum_{j>N}\|\Phi g_j\|_H^2\right)\mathrm{d}t\right]^{1/2}
\longrightarrow0, ~\text{as}~ N\rightarrow\infty.
\end{aligned}
\]
Hence,
\[\sigma(u_N)-\bar{\sigma}_N(u_N)\rightharpoonup 0 ~\text{weakly in }~L^2\bigl(\Omega;L^2(0,T;\mathcal{L}_2(\mathcal{U},H))\bigr).\]
Combining this with \eqref{3.17}, we conclude that
\begin{align}\label{3.17jia}
\sigma(u_N)\rightharpoonup \chi_1~\text{weakly in } L^2\bigl(\Omega;L^2(0,T;\mathcal{L}_2(\mathcal{U},H))\bigr).
\end{align}
Moreover, since $\int_0^t P_N\sigma(u_N(r))\mathrm{d}W_N(r)=\int_0^t \bar{\sigma}_N(u_N(r))\mathrm{d}W(r)$
and the stochastic integral defines a bounded linear mapping from $L^2\bigl(\Omega;L^2(0,T;L_2(\mathcal{U},H))\bigr)$
into $L^2(\Omega;H)$, it follows from \eqref{3.17} that
\begin{align}\label{3.18}
\int_0^t P_N\sigma(u_N(r))\mathrm{d}W_N(r)\rightharpoonup\int_0^t \chi_1(r)\mathrm{d}W(r)\quad\text{weakly in }L^2(\Omega;H).
\end{align}
Next, we consider the limits of each term in \eqref{lili1} as $N\rightarrow\infty$ one by one.
Taking into account $\eqref{lili1}_1$, together with integration by parts and Fubini's theorem, it follows that for every $\varrho\in L^{\infty}([0,T]\times \Omega,\mathbb{R})$, $j\in\mathbb N$  and $v=e_j$, for $N\geq j$, we have
\begin{align}\label{zhuo}
&\mathbb{E}\int_{0}^{T}\varrho(t)\langle u_{N}(t),e_{j}\rangle_{H}\mathrm{d}t-\mathbb{E}\int_{0}^{T}\varrho(t)\langle P_{N}u_{0},e_{j}\rangle_{H}\mathrm{d}t\notag\\
=&\mathbb{E}\left(\int_{0}^{T}\varrho(t)\int_{0}^{t}\left(-\kappa a_{\alpha}( u_{N}(r),e_{j})-(1-\kappa)\int_{0}^{\infty}\mu(s)a_{\alpha} ( \eta_{N}(r,s),e_{j})\mathrm{d}s\notag\right)\mathrm{d}r\mathrm{d}t\right)\notag\\
&-\lambda\mathbb{E}\int_{0}^{T}\varrho(t)\int_{0}^{t}\big\langle u_{N}(r),e_{j}\big\rangle_{H}\mathrm{d}r\mathrm{d}t+\mathbb{E}\int_{0}^{T}\varrho(t)\int_{0}^{t}\big\langle P_{N}\varphi(u_{N}(r)),e_{j}\big\rangle_{\langle L^{q}(\mathcal{O}),L^{p_{0}+1}(\mathcal{O})\rangle}\mathrm{d}r\mathrm{d}t\notag\\
&+\mathbb{E}\int_{0}^{T}\varrho(t)\int_{0}^{t}\big\langle P_{N}\sigma(u_{N}(r))\mathrm{d}W_{N}(r),e_{j}\big\rangle_{H}\mathrm{d}t\notag\\
=&\mathbb{E}\int_{0}^{T}-\kappa a_{\alpha}(u_{N}(t),e_{j})\int_{t}^{T}\varrho(r)\mathrm{d}r\mathrm{d}t\notag\\
&+\mathbb{E}\int_{0}^{T} \Big\langle P_{N}\varphi(u_{N}(t)),e_{j}\int_{t}^{T}\varrho(r)\mathrm{d}r\Big\rangle_{\langle L^{q}(\mathcal{O}),L^{p_{0}+1}(\mathcal{O})\rangle}\mathrm{d}t\notag\\
&-(1-\kappa)\mathbb{E}\int_{0}^{T}\int_{0}^{\infty}\mu(s)a_{\alpha} (\eta_{N}(t,s),e_{j}) \left(\int_{t}^{T}\varrho(r)\mathrm{d}r\right) \mathrm{d}s\mathrm{d}t\notag\\
&-\lambda\mathbb{E}\int_{0}^{T} \Big\langle u_{N}(t),e_{j}\int_{t}^{T}\varrho(r)\mathrm{d}r\Big\rangle_{H}\mathrm{d}t+\int_{0}^{T}\mathbb{E}\Big\langle\varrho(t)e_{j},\int_{0}^{t}P_{N}\sigma(u_{N}(r))\mathrm{d}W_{N}(r)\Big\rangle_{H}\mathrm{d}t.
\end{align}
Letting $N\rightarrow\infty$ in \eqref{zhuo}, by \eqref{zhuo0}-\eqref{3.18}, we have
\begin{align}\label{zhuo1}
&\mathbb{E}\int_{0}^{T}\varrho(t)\langle u(t),e_{j}\rangle_{H}\mathrm{d}t-\mathbb{E}\int_{0}^{T}\varrho(t)\langle u_{0},e_{j}\rangle_{H}\mathrm{d}t\notag\\
=&\mathbb{E}\left(\int_{0}^{T}\varrho(t)\int_{0}^{t}\left(-\kappa a_{\alpha}( u(r),e_{j})-(1-\kappa)\int_{0}^{\infty}\mu(s)a_{\alpha}( \eta(r,s),e_{j})\mathrm{d}s\notag\right)\mathrm{d}r\mathrm{d}t\right)\notag\\
&+\mathbb{E}\int_{0}^{T}\varrho(t)\int_{0}^{t}\langle \chi(r),e_{j}\rangle_{\langle L^{q}(\mathcal{O}),L^{p_{0}+1}(\mathcal{O})\rangle}\mathrm{d}r\mathrm{d}t-\lambda\mathbb{E}\int_{0}^{T}\varrho(t)\int_{0}^{t}\langle u(r),e_{j}\rangle_{H}\mathrm{d}r\mathrm{d}t\notag\\
&+\mathbb{E}\int_{0}^{T}\varrho(t)\int_{0}^{t}\langle\chi_{1}(r)\mathrm{d}W(r),e_{j}\rangle_{H}\mathrm{d}t.
\end{align}
Since $\operatorname{span}\{e_j:j\in\mathbb N\}$ is dense in $V\cap L^{p_0+1}(\mathcal O)$, similar to \cite{NEW1}, we can extend \eqref{zhuo1}  to every $v\in V\cap L^{p_0+1}(\mathcal O)$.
Immediately, we find that
\begin{align}\label{zhuo2}
u(t)=&u_{0}-\kappa\int_{0}^{t}A_{\alpha}u(r)\mathrm{d}r-\lambda\int_{0}^{t}u(r)\mathrm{d}r-(1-\kappa)\int_{0}^{t}\int_{0}^{\infty}\mu(s)A_{\alpha}\eta(r,s)\mathrm{d}s\mathrm{d}r\notag\\
&+\int_{0}^{t}\chi(r)\mathrm{d}r+\int_{0}^{t}\chi_{1}(r)\mathrm{d}W(r),
\end{align}
in $(V\cap L^{p_{0}+1}(\mathcal{O}))^{*},$ $\mathrm{d}t\times \mathbb{P}$-a.s. on $[0,T]\times\Omega$.

Next, we identify the limit of $\eqref{lili1}_{2}$ as $N\rightarrow\infty$. For  equation  $\eqref{lili1}_{2}$,
 by expression \eqref{x}, for all $\tilde{\eta} \in M_{\mu}^{0}$, we find that
\begin{align*}
\langle \eta_N(t),\tilde{\eta}\rangle_{M_{\mu}^{0}}
=&\int_{0}^{t} \int_{0}^{s}
a_{\alpha}(u_N(t-r),\tilde{\eta}(s))\mathrm{d}r\mu(s)\mathrm{d}s
+\int_{t}^{\infty}
a_{\alpha}( R_N\eta_0(s-t),\tilde{\eta}(s))\mu(s)\mathrm{d}s\\
&+\int_{t}^{\infty}\int_{0}^{t}a_{\alpha}(u_N(t-r),\tilde{\eta}(s))\mathrm{d}r\mu(s)\mathrm{d}s.
\end{align*}
Taking limit as $N \to \infty$, we arrive at
\begin{align}\label{shui}
\langle\eta(t),\tilde{\eta}\rangle_{M_{\mu}^{0}}=&\int_{0}^{t}\int_{0}^{s}a_{\alpha}(u(t-r),\tilde{\eta}(s))\mathrm{d}r\mu(s)\mathrm{d}s
+\int_{t}^{\infty}a_{\alpha}(\eta_0(s-t),\tilde{\eta}(s))\mu(s)\mathrm{d}s\nonumber\\
&+\int_{t}^{\infty}\int_{0}^{t}a_{\alpha}(u(t-r),\tilde{\eta}(s))\mathrm{d}r\mu(s)\mathrm{d}s.
\end{align}
Since the above identity holds for all
$\tilde{\eta}$, $\eta(t)$ indeed satisfies expression \eqref{x}, which is equivalent to equation \eqref{nan}, i.e.,
\begin{align}\label{qa}
\partial_t\eta(t)=\mathcal{T}_{\mu}\eta(t)+u(t).
\end{align}
in the weak sense in $M_\mu^0$.

 In what follows, we prove $\chi=\varphi(u)$ and $\chi_{1}=\sigma(u)$, $\mathrm{d}t\times \mathbb{P}$-a.s. on $[0,T]\times\Omega$. Now, let $\hat{u}=u-\varepsilon\hat{\varrho}\xi$,  where $\varepsilon\in[0,1]$, $\hat{\varrho}(0)=0$, $\hat{\varrho}\in L^{\infty}([0,T]\times\Omega,\mathbb{R})$ and $\xi\in H^{\alpha}\cap L^{p_{0}+1}(\mathcal{O})$.
Define $\hat{\eta}$ as the solution of the linear equation
\[
\partial_t \hat{\eta}=\mathcal{T}_\mu \hat{\eta}+\hat{u},
\qquad
\hat{\eta}(0,s)=\eta_0(s) .
\]
%Then, we have
%\begin{align}\label{lyu}
%\partial_t(\eta-\hat\eta)=\mathcal{T}_\mu(\eta-\hat\eta)+(u-\hat u),\qquad(\eta-\hat\eta)(0,s)=0.
%\end{align}
Define $G=2a_{\varphi}^{+}+L_\sigma^2+1$, where $\alpha_{\varphi}^{+}=\max\{0,\alpha_{\varphi}\}$. Then, applying  It\^{o}'s formula,
%%%to
%%%\textcolor{red}{CHECK}
%%%\[
%%%\begin{aligned}
%%% e^{-Gt}\mathrm{d}u_{N}(t)=&e^{-Gt}\Bigl(-\kappa(-\Delta)^{\alpha}u_{N}(t)\mathrm{d}t-\lambda u_{N}(t)\mathrm{d}t
%%%-(1-\kappa)\int_0^\infty \mu(s)(-\Delta)^{\alpha}\eta_{N}(t,s)\mathrm{d}s\mathrm{d}t\\
%%%&+P_N\varphi(u_{N}(t))\mathrm{d}t
%%%+\sigma_N(u_{N}(t))\mathrm{d}W(t)\Bigr),
%%%\end{aligned}
%%%\]
we have
\begin{align}\label{zhuo5}
&e^{-Gt}\mathbb{E}\|u_{N}(t)\|_{H}^{2}-\|P_{N}u_{0}\|_{H}^{2}\notag\\[1.0ex]
=&\mathbb{E}\Bigg[\int_{0}^{t} e^{-Gs}
\Big(-2\kappa a_{\alpha}(u_{N}(s),u_{N}(s))
-2(1-\kappa)\int_{0}^{\infty}\mu(k)a_{\alpha}(\eta_{N}(s,k),u_{N}(s))\mathrm{d}k\notag\\[1.0ex]
&-2\lambda\langle u_{N}(s),u_{N}(s)\rangle_{H}+2\big\langle P_{N}\varphi(u_{N}(s)),u_{N}(s)\big\rangle_{\langle L^{q}(\mathcal{O}),L^{p_{0}+1}(\mathcal{O})\rangle}\notag\\[0.8ex]
&+\|\bar{\sigma}_{N}(u_{N}(s))\|_{\mathcal{L}_{2}(\mathcal{U},H)}^{2}-G\|u_{N}(s)\|_{H}^{2}\Big)\mathrm{d}s\Bigg].
\end{align}
Since $P_N$ is an orthogonal projection and $\bar{\sigma}_N(u_N)$ is the zero extension of $P_N\sigma(u_N):\mathcal{U}_N\to H_N$, we have
\[\|\bar{\sigma}_N(u_N)\|_{\mathcal{L}_2(\mathcal{U},H)}^2=\sum_{j=1}^{N}\|P_N\sigma(u_N)g_j\|_H^2\leq\|\sigma(u_N)\|_{\mathcal{L}_2(\mathcal{U},H)}^2.\]
Therefore, we have
\begin{align}\label{zhuoq}
&e^{-Gt}\mathbb{E}\|u_{N}(t)\|_{H}^{2}-\|P_{N}u_{0}\|_{H}^{2}\notag\\[1.0ex]
\leq&\mathbb{E}\Bigg[\int_{0}^{t} e^{-Gs}
\Big(-2\kappa a_{\alpha}(u_{N}(s),u_{N}(s))
-2(1-\kappa)\int_{0}^{\infty}\mu(k)a_{\alpha}(\eta_{N}(s,k),u_{N}(s))\mathrm{d}k\notag\\[1.0ex]
&-2\lambda\langle u_{N}(s),u_{N}(s)\rangle_{H}+2\big\langle \varphi(u_{N}(s)),u_{N}(s)\big\rangle_{\langle L^{q}(\mathcal{O}),L^{p_{0}+1}(\mathcal{O})\rangle}\notag\\[0.8ex]
&+\|\sigma(u_{N}(s))\|_{\mathcal{L}_{2}(\mathcal{U},H)}^{2}-G\|u_{N}(s)\|_{H}^{2}\Big)\mathrm{d}s\Bigg].
\end{align}
Inserting $\hat{u}$ and $\hat{\eta}$ into the above inequality, we obtain
\begin{align}\label{zhuo500}
&e^{-Gt}\mathbb{E}\|u_{N}(t)\|_{H}^{2}-\|P_{N}u_{0}\|_{H}^{2}\notag\\[1.2ex]
\leq&\mathbb{E}\bigg[\int_{0}^{t}e^{-Gs}\bigg(-2\kappa a_{\alpha}(\hat{u}(s),u_{N}(s))
-2\kappa a_{\alpha}((u_{N}(s)-\hat{u}(s)),\hat{u}(s))\notag\\
&\left.-2(1-\kappa)\int_{0}^{\infty}\mu(k)a_{\alpha}(\hat{\eta}(s,k),u_{N}(s))\mathrm{d}k\right.\notag\\
&\left.-2\lambda\langle\hat{u}(s),u_{N}(s)\rangle_{H}+2\lambda\langle-u_{N}(s)+\hat{u}(s),\hat{u}(s)\rangle_{H}\right.\notag\\[1.4ex]
&\left.+2(1-\kappa)\int_{0}^{\infty}\mu(k)a_{\alpha}( - \eta_{N}(s,k)+\hat{\eta}(s,k),\hat{u}(s))\mathrm{d}k\right.\notag\\[1.2ex]
&\left.+2\langle \varphi(\hat{u}(s)),u_{N}(s)\rangle_{\langle L^{q}(\mathcal{O}), L^{p_{0}+1}(\mathcal{O})\rangle}+2\langle \varphi (u_{N}(s))-\varphi(\hat{u}(s)),\hat{u}(s)\rangle_{\langle L^{q}(\mathcal{O}), L^{p_{0}+1}(\mathcal{O})\rangle}\right.\notag\\[1.2ex]
&\left.-\|\sigma(\hat{u}(s))\|_{\mathcal{L}_{2}(\mathcal{U},H)}^{2}+2\langle \sigma(u_{N}(s)),\sigma(\hat{u}(s))\rangle_{\mathcal{L}_{2}(\mathcal{U},H)}\right.\notag\\[1.2ex]
&-G\left(2\langle u_{N}(s),\hat{u}(s)\rangle_{H}-\|\hat{u}(s)\|_{H}^{2}\right)
\bigg)\mathrm{d}s\bigg]\notag\\
&+\mathbb{E}\bigg[\int_{0}^{t}e^{-Gs}\times\bigg(-2\kappa a_{\alpha}((u_{N}(s)-\hat{u}(s)),(u_{N}(s)-\hat{u}(s))) \notag\\[1.2ex]
&\left.+2(1-\kappa)\underbrace{\int_{0}^{\infty}\mu(k)a_{\alpha}(- \eta_{N}(s,k)+\hat{\eta}(s,k),(u_{N}(s)-\hat{u}(s)))\mathrm{d}k}
_{=D}\right.\notag\\[1.2ex]
&\left.+2\langle \varphi(u_{N}(s))-\varphi(\hat{u}(s)),u_{N}(s)-\hat{u}(s)\rangle_{\langle L^{q}(\mathcal{O}), L^{p_{0}+1}(\mathcal{O})\rangle}-G\|u_{N}(s)-\hat{u}(s)\|_{H}^{2}\right.\notag\\[1.2ex]
&+\|\sigma(u_{N}(s))-\sigma(\hat{u}(s))\|_{\mathcal{L}_{2}(\mathcal{U},H)}^{2}-2\lambda\|u_{N}(s)-\hat{u}(s)\|_{H}^{2}\bigg)\mathrm{d}s\bigg].
\end{align}
We now proceed to the discussion of $D$. Let $w:=u_{N}-\hat u,~\zeta:=\eta_{N}-\hat\eta.$
By Lemma \ref{mon}, we have
\[\partial_t\zeta=\mathcal{T}_\mu\zeta+w,\qquad\zeta(0,s)=R_{N}\eta_{0}(s)-\eta_{0}(s).\]
Taking the inner product of the above equation with \(\zeta\) in \(M_\mu^0\), we obtain
\[\frac{1}{2}\frac{\mathrm{d}}{\mathrm{d}t}\|\zeta\|_{M_\mu^0}^2=\langle \mathcal{T}_\mu\zeta,\zeta\rangle_{M_\mu^0}+\langle w,\zeta\rangle_{M_\mu^0}.\]
Hence,
\begin{align*}
\mathbb{E}\int_0^t e^{-Gs}D\mathrm{d}s=&-2(1-\kappa)\mathbb{E}\int_0^te^{-Gs}\langle w(s),\zeta(s)\rangle_{M_\mu^0}\mathrm{d}s\\
=&\mathbb{E}\int_0^t e^{-Gs}\left[-(1-\kappa)\frac{d}{ds}\|\zeta(s)\|_{M_\mu^0}^2+2(1-\kappa)\langle T_\mu\zeta(s),\zeta(s)\rangle_{M_\mu^0}\right]\mathrm{d}s\\
\leq&(1-\kappa)\mathbb{E}\Bigg[-e^{-Gt}\|\zeta(t)\|_{M_\mu^0}^2+\|\zeta(0)\|_{M_\mu^0}^2\\
&-G\int_0^t e^{-Gs}\|\zeta(s)\|_{M_\mu^0}^2\mathrm{d}s-\delta\int_0^t e^{-Gs}\|\zeta(s)\|_{M_\mu^0}^2\mathrm{d}s\Bigg].
\end{align*}
It follows from the above  expression, $(P_{3})$ and $(Q_{0})$ that
\begin{align}\label{zhuo500}
&e^{-Gt}\mathbb{E}\|u_{N}(t)\|_{H}^{2}-\|P_{N}u_{0}\|_{H}^{2}\notag\\[1.2ex]
\leq&\mathbb{E}\bigg[\int_{0}^{t}e^{-Gs}\bigg(-2\kappa a_{\alpha}(\hat{u}(s),u_{N}(s))
-2\kappa a_{\alpha}((u_{N}(s)-\hat{u}(s)),\hat{u}(s))\notag\\
&\left.-2(1-\kappa)\int_{0}^{\infty}\mu(k)a_{\alpha}(\hat{\eta}(s,k),u_{N}(s))\mathrm{d}k\right.\notag\\
&\left.-2\lambda\langle\hat{u}(s),u_{N}(s)\rangle_{H}+2\lambda\langle-u_{N}(s)+\hat{u}(s),\hat{u}(s)\rangle_{H}\right.\notag\\[1.4ex]
&\left.+2(1-\kappa)\int_{0}^{\infty}\mu(k)a_{\alpha}(-\eta_{N}(s,k)+\hat{\eta}(s,k),\hat{u}(s))\mathrm{d}k\right.\notag\\[1.2ex]
&\left.+2\langle \varphi(\hat{u}(s)),u_{N}(s)\rangle_{\langle L^{q}(\mathcal{O}), L^{p_{0}+1}(\mathcal{O})\rangle}+2\langle \varphi (u_{N}(s))-\varphi(\hat{u}(s)),\hat{u}(s)\rangle_{\langle L^{q}(\mathcal{O}), L^{p_{0}+1}(\mathcal{O})\rangle}\right.\notag\\[1.2ex]
&\left.-\|\sigma(\hat{u}(s))\|_{\mathcal{L}_{2}(\mathcal{U},H)}^{2}+2\langle \sigma(u_{N}(s)),\sigma(\hat{u}(s))\rangle_{\mathcal{L}_{2}(\mathcal{U},H)}\right.\notag\\[1.2ex]
&-G\left(2\langle u_{N}(s),\hat{u}(s)\rangle_{H}-\|\hat{u}(s)\|_{H}^{2}\right)
\bigg)\mathrm{d}s\bigg]+(1-\kappa)\|\zeta(0)\|_{M_\mu^0}^2.
\end{align}
Therefore, letting $N\rightarrow\infty$, we see that for any $\varepsilon_{1}>0$, there exists a positive constant $\tilde{N}$ such that, for all $N>\tilde{N}$,
$$(1-\kappa)\|\zeta(0)\|_{M_\mu^0}^2\leq\varepsilon_{1}.$$
Moreover, it follows that for any nonnegative $\varrho\in L^{\infty}([0,T],\mathbb{R})$,
\begin{align}\label{chao}
&\int_{0}^{T}\varrho(t)\left(e^{-Gt}\mathbb{E}\|u(t)\|_{H}^{2}-\|u_{0}\|_{H}^{2}\right)\mathrm{d}t\notag\\
\leq&\liminf_{N\rightarrow\infty}\int_{0}^{T}\varrho(t)\left(e^{-Gt}\mathbb{E}\|u_{N}(t)\|_{H}^{2}-\|P_{N}u_{0}\|_{H}^{2}\right)\mathrm{d}t\notag\\[0.8ex]
\leq&\mathbb{E}\bigg[\int_{0}^{T}\varrho(t)\int_{0}^{t}e^{-Gs}\bigg(-2\kappa a_{\alpha}(\hat{u}(s),u(s))
-2\kappa a_{\alpha}((u(s)\notag\\
&\left.-\hat{u}(s)),\hat{u}(s))-2(1-\kappa)\int_{0}^{\infty}\mu(k)a_{\alpha}(\hat{\eta}(s,k),u(s))\mathrm{d}k\right.\notag\\
&\left.
+2(1-\kappa)\int_{0}^{\infty}\mu(k)a_{\alpha}(- \eta(s,k)+\hat{\eta}(s,k),\hat{u}(s))\mathrm{d}k\right.\notag\\[0.8ex]
&\left.+2\lambda\langle-u(s)+\hat{u}(s),\hat{u}(s)\rangle_{H}-2\lambda\langle\hat{u}(s),u(s)\rangle_{H}
+2\langle \varphi(\hat{u}(s)),u(s)\rangle_{\langle L^{q}(\mathcal{O}), L^{p_{0}+1}(\mathcal{O})\rangle}\right.\notag\\[1.4ex]
&\left.+2\langle \chi-\varphi(\hat{u}(s)),\hat{u}(s)\rangle_{\langle L^{q}(\mathcal{O}), L^{p_{0}+1}(\mathcal{O})\rangle}-\|\sigma(\hat{u}(s))\|_{\mathcal{L}_{2}(\mathcal{U},H)}^{2}
+2\langle \chi_{1},\sigma(\hat{u}(s))\rangle_{\mathcal{L}_{2}(\mathcal{U},H)}\right.\notag\\[1.0ex]
&-G\left(2\langle u(s),\hat{u}(s)\rangle_{H}-\|\hat{u}(s)\|_{H}^{2}\right)
\bigg)\mathrm{d}s\mathrm{d}t\bigg]+\varepsilon_{1}\int_{0}^{T}\varrho(t)\mathrm{d}t.
\end{align}
On the other hand, using It\^{o}'s formula to \eqref{zhuo2}, we find
\begin{equation}\label{zhuo7}
\begin{aligned}
&\int_{0}^{T}\varrho(t)\left(e^{-Gt}\mathbb{E}\|u(t)\|_{H}^{2}-\|u_{0}\|_{H}^{2}\right)\mathrm{d}t\\
=&\mathbb{E}\Bigl(\int_{0}^{T}\varrho(t)\int_{0}^{t}e^{-Gs}\Bigl(-2\kappa a_{\alpha}(u(s),u(s))\\
&-2(1-\kappa)\int_{0}^{\infty}\mu(k)a_{\alpha}(\eta(s,k),u(s))\mathrm{d}k\mathrm{d}s\mathrm{d}t\Bigr)\\
&+\mathbb{E}\Bigl(\int_{0}^{T}\varrho(t)\int_{0}^{t}e^{-Gs}\Bigl(-2\lambda\|u(s)\|_{H}^{2}+2\langle\chi,u(s)\rangle_{\langle L^{q}(\mathcal{O}), L^{p_{0}+1}(\mathcal{O})\rangle}\\
&+\|\chi_{1}\|_{\mathcal{L}_{2}(\mathcal{U},H)}^{2}-G\|u(s)\|_{H}^{2}\Bigr)\mathrm{d}s\mathrm{d}t\Bigr).
\end{aligned}
\end{equation}
Substituting the above equation into the left-hand side of \eqref{chao}, rearranging the resulting inequality, and then letting $\varepsilon_1\to 0$, we obtain
\begin{align}\label{zhuo8}
0\geq&\mathbb{E}\bigg[\int_{0}^{T}\varrho(t)\int_{0}^{t}e^{-Gs}\left(
-2\kappa a_{\alpha}((u(s)-\hat{u}(s)),(u(s)-\hat{u}(s)))\right.\notag\\
&\left.+2(1-\kappa)\int_{0}^{\infty}\mu(k)a_{\alpha}( - \eta(s,k)+ \hat{\eta}(s,k),(u(s)-\hat{u}(s)))\mathrm{d}k\right.\notag\\
&\left.+2\langle\chi-\varphi(\hat{u}(s)),u(s)-\hat{u}(s)\rangle_{\langle L^{q}(\mathcal{O}), L^{p_{0}+1}(\mathcal{O})\rangle}-2\lambda \| u(s)-\hat{u}(s)\|_{H}^{2}\right.\notag\\[0.8ex]
&\left.+\|\chi_{1}-\sigma(\hat{u}(s))\|_{\mathcal{L}_{2}(\mathcal{U},H)}^{2}-G\|u(s)-\hat{u}(s)\|_{H}^{2}
\right)\mathrm{d}s\mathrm{d}t\bigg].
\end{align}
Recall that $\hat{u}=u-\varepsilon\hat{\varrho}\xi,$ where $\varepsilon\in[0,1],~\hat{\varrho}\in L^\infty([0,T]\times\Omega;\mathbb{R}),~
\xi\in V\cap L^{p_0+1}(\mathcal O).$
In particular, when $\varepsilon=0$, we have $\hat{u}=u$. Hence, by the uniqueness of the solution to the memory equation, we also have
$\hat{\eta}=\eta.$
Therefore, it follows from \eqref{zhuo8} that
\[\chi_1=\sigma(u),~\mathrm{d}t\times\mathbb{P}\text{-a.s. on }[0,T]\times\Omega.\]
Now, we prove $\chi=\varphi(u)$, $\mathrm{d}t\times\mathbb P$-a.s. on
$[0,T]\times\Omega$. It follows from
$\hat u=u-\varepsilon\hat\varrho\xi$
 that
\begin{align}\label{to}
0\geq&\mathbb{E}\bigg[\int_{0}^{T}\varrho(t)\int_{0}^{t}e^{-Gs}\varepsilon\hat{\varrho}(s)\left(
2\kappa a_{\alpha}(-u(s)+\hat{u}(s),\xi)\right.\notag\\
&\left.+2(1-\kappa)\int_{0}^{\infty}\mu(k)a_{\alpha}(-\eta(s,k)+\hat{\eta}(s,k),\xi)\mathrm{d}k\right.\notag\\
&\left.+2\langle\chi-\varphi(\hat{u}(s)),\xi\rangle_{\langle L^{q}(\mathcal{O}), L^{p_{0}+1}(\mathcal{O})\rangle}
\right)\mathrm{d}s\mathrm{d}t\bigg]
+\mathbb{E}\bigg[\int_{0}^{T}\varrho(t)\int_{0}^{t}e^{-Gs}\left(-2\lambda \varepsilon^{2}(\hat{\varrho}(s))^{2}\|\xi\|_{H}^{2}\right.\notag\\[0.8ex]
&\left.+\|\chi_{1}-\sigma(\hat{u}(s))\|_{\mathcal{L}_{2}(\mathcal{U},H)}^{2}-G\varepsilon^{2}(\hat{\varrho}(s))^{2}\|\xi\|_{H}^{2}
\right)\mathrm{d}s\mathrm{d}t\bigg].
\end{align}
After dividing the above inequality by $\varepsilon>0$, we obtain
\begin{align}\label{3.32}
0\geq&\mathbb{E}\bigg[\int_{0}^{T}\varrho(t)\int_{0}^{t}e^{-Gs}\hat{\varrho}(s)\left(
2\kappa a_{\alpha}(-u(s)+\hat{u}(s),\xi)\right.\notag\\
&\left.+2(1-\kappa)\int_{0}^{\infty}\mu(k)a_{\alpha}(-\eta(s,k)+\hat{\eta}(s,k),\xi)\mathrm{d}k\right.\notag\\
&\left.+2\langle\chi-\varphi(\hat{u}(s)),\xi\rangle_{\langle L^{q}(\mathcal{O}), L^{p_{0}+1}(\mathcal{O})\rangle}
\right)\mathrm{d}s\mathrm{d}t\bigg]
+\mathbb{E}\bigg[\int_{0}^{T}\varrho(t)\int_{0}^{t}e^{-Gs}\left(-2\lambda \varepsilon(\hat{\varrho}(s))^{2}\|\xi\|_{H}^{2}\right.\notag\\[0.8ex]
&\left.+\frac{1}{\varepsilon}\|\chi_{1}-\sigma(\hat{u}(s))\|_{\mathcal{L}_{2}(\mathcal{U},H)}^{2}-G\varepsilon(\hat{\varrho}(s))^{2}\|\xi\|_{H}^{2}
\right)\mathrm{d}s\mathrm{d}t\bigg].
\end{align}
Define $F_\varepsilon(s):=2(1-\kappa)\int_0^\infty \mu(k)a_{\alpha}(-\eta(s,k)+\hat{\eta}(s,k),\xi)\mathrm{d}k$
and \\
$R_\varepsilon:=\mathbb{E}\Bigg[\int_0^T \varrho(t)\int_0^t e^{-Gs}\Bigg(2\kappa \hat{\varrho}(s)a_{\alpha}(-u(s)
+\hat{u}(s),\xi)-2\lambda\varepsilon\bigl(\hat{\varrho}(s)\bigr)^2\|\xi\|_H^2
+\frac{1}{\varepsilon}\|\sigma(u(s))-\sigma(\hat{u}(s))\|_{\mathcal{L}_2(U,H)}^2-G\varepsilon\bigl(\hat{\varrho}(s)\bigr)^2\|\xi\|_H^2
\Bigg)\mathrm{d}s\mathrm{d}t\Bigg]$.

Therefore, \eqref{3.32} can be rewritten as
\begin{align}\label{3.33}
0 \geq \mathbb{E}\Bigg[\int_0^T \varrho(t)\int_0^t e^{-Gs}\hat{\varrho}(s)\Bigg(F_\varepsilon(s)+2\left\langle\chi(s)-\varphi(\hat{u}(s)),\xi\right\rangle_{L^q(\mathcal O),L^{p_0+1}(\mathcal O)}\Bigg)\mathrm{d}s\mathrm{d}t\Bigg]+R_\varepsilon .
\end{align}
Since $\hat{u}(s)-u(s)=-\varepsilon\hat{\varrho}(s)\xi,$ we have
\[2\kappa\hat{\varrho}(s)a_{\alpha}(-u(s)+\hat{u}(s),\xi)=
-2\kappa\varepsilon\bigl(\hat{\varrho}(s)\bigr)^2a_{\alpha}(\xi,\xi).\]
Moreover, by $(Q_0)$, we have
\[\frac{1}{\varepsilon}\left\|\sigma(u(s))-\sigma(\hat{u}(s))\right\|_{\mathcal{L}_2(\mathcal{U},H)}^2\leq\frac{L_\sigma^2}{\varepsilon}\|u(s)-\hat{u}(s)\|_H^2
=L_\sigma^2\varepsilon\bigl(\hat{\varrho}(s)\bigr)^2\|\xi\|_H^2.\]
Therefore,
\[R_\varepsilon \longrightarrow 0~\text{as }\varepsilon\to0,\]
and
\begin{align}\label{3.33}
0 \geq \mathbb{E}\Bigg[\int_0^T \varrho(t)\int_0^t e^{-Gs}\hat{\varrho}(s)\Bigg(F_\varepsilon(s)+2\Big\langle\chi(s)-\varphi(\hat{u}(s)),\xi
\Big\rangle_{L^q(\mathcal O),L^{p_0+1}(\mathcal O)}\Bigg)\mathrm{d}s\mathrm{d}t\Bigg].
\end{align}
By Lemma \ref{mon}, we have
\[\eta(t,s)-\hat{\eta}(t,s)
=\begin{cases}
\displaystyle
\int_0^s \bigl(u(t-r)-\hat{u}(t-r)\bigr)\mathrm{d}r,
& 0<s\leq t, \\[1.2ex]
\displaystyle
\int_0^t \bigl(u(t-r)-\hat{u}(t-r)\bigr)\mathrm{d}r,
& s>t.
\end{cases}
\]
Moreover, since
$u-\hat{u}=\varepsilon\hat{\rho}\xi,$
we obtain
\[\eta(t,s)-\hat{\eta}(t,s)
=\varepsilon
\begin{cases}
\displaystyle
\int_0^s \hat{\rho}(t-r)\xi\mathrm{d}r,
& 0<s\leq t, \\[1.2ex]
\displaystyle
\int_0^t \hat{\rho}(t-r)\xi\mathrm{d}r,
& s>t.
\end{cases}
\]
Hence, $F_\varepsilon(s)\to0~\text{as }\varepsilon\to0.$
Using the continuity and growth condition of $\varphi$, together with the dominated convergence theorem, and letting $\varepsilon\to0$ in \eqref{3.33} we obtain
\begin{align}\label{3.34}
0 \geq2\mathbb{E}\Bigg[\int_0^T \varrho(t)\int_0^t e^{-Gs}\hat{\varrho}(s)\left\langle\chi(s)-\varphi(u(s)),\xi\right\rangle_{L^q(\mathcal O),L^{p_0+1}(\mathcal O)}
\mathrm{d}s\mathrm{d}t\Bigg].
\end{align}
By the arbitrariness of $\varrho$, $\hat{\varrho}$, and $\xi$, we conclude that
\[\chi=\varphi(u)~\text{in }L^q(\mathcal O),~\mathrm{d}t\times\mathbb{P}\text{-a.s. on }[0,T]\times\Omega.\]
This, together with \eqref{zhuo2} and \eqref{qa}, implies that $(u,\eta)$ is a solution to \eqref{lii10} in the sense of Definition \ref{dingyi1}.

\textbf{Step 4. Uniqueness and continuous with respect to initial data.} Let $(u_1,\eta_1)$ and $(u_2,\eta_2)$ be two solutions of problem \eqref{lii10} corresponding to the same initial data $U_0=(u_0,\eta_0).$ Set $z:=u_1-u_2,~\zeta:=\eta_1-\eta_2.$
Then $(z,\zeta)$ satisfies
\[\left\{
\begin{aligned}\label{3.35}
&\mathrm{d}z(t)+\kappa A_\alpha z(t)\mathrm{d}t+(1-\kappa)\int_0^\infty\mu(s)A_\alpha\zeta(t,s)\mathrm{d}s\mathrm{d}t
-\bigl(\varphi(u_1(t))-\varphi(u_2(t))\bigr)\mathrm{d}t+\lambda z(t)\mathrm{d}t\\
&=\bigl(\sigma(u_1(t))-\sigma(u_2(t))\bigr)\mathrm{d}W(t),\\
&\partial_t\zeta(t)=\mathcal{T}_\mu\zeta(t)+z(t),\\
&z(0)=0,~\zeta(0,s)=0,~\zeta(t,0)=0.
\end{aligned}
\right.
\]
Applying It\^{o}'s formula to $\Psi_0(z,\zeta)=\frac12\|z\|_H^2+\frac{1-\kappa}{2}\|\zeta\|_{M_\mu^0}^2,$ we obtain
\[\begin{aligned}
\mathrm{d}\Psi_0(z(t),\zeta(t))
=&-\kappa a_{\alpha}(z(t),z(t))\mathrm{d}t-\lambda\|z(t)\|_H^2\mathrm{d}t+
(1-\kappa)\langle \mathcal{T}_\mu\zeta(t),\zeta(t)\rangle_{M_\mu^0}\mathrm{d}t\\
&+\left\langle\varphi(u_1(t))-\varphi(u_2(t)),z(t)\right\rangle_{L^q(\mathcal O),L^{p_0+1}(\mathcal O)}\mathrm{d}t
+\frac12\|\sigma(u_1(t))-\sigma(u_2(t))\|_{L_2(\mathcal{U},H)}^2\mathrm{d}t\\
&+\left\langle z(t),\bigl(\sigma(u_1(t))-\sigma(u_2(t))\bigr)\mathrm{d}W(t)\right\rangle_H .
\end{aligned}
\]
Using the dissipativity of $\mathcal{T}_\mu$, we have
\[\langle \mathcal{T}_\mu\zeta,\zeta\rangle_{M_\mu^0}\leq-\frac{\delta}{2}\|\zeta\|_{M_\mu^0}^2,\]
together with assumptions $(P_3)$ and $(Q_0)$, we obtain
\[\mathbb E\Psi_0(z(t),\zeta(t))\leq C\int_0^t\mathbb E\Psi_0(z(s),\zeta(s))\mathrm{d}s.\]
Since $\Psi_0(z(0),\zeta(0))=0,$ by Gronwall's inequality, we know
\[\mathbb E\Psi_0(z(t),\zeta(t))=0,~t\in[0,T].\]
Hence $z=0,~\zeta=0,~\mathrm{d}t\times\mathbb P\text{-a.s. on }[0,T]\times\Omega,$ therefore the solution is unique.

To prove continuous dependence on initial data, let $U_0=(u_0,\eta_0),~\widetilde U_0=(\widetilde u_0,\widetilde\eta_0),$
and let $U(t)=(u(t),\eta(t)),~\widetilde U(t)=(\widetilde u(t),\widetilde\eta(t))$
be the corresponding solutions. Repeating the above argument for the
difference of the two solutions and using the BDG inequality together
with Gronwall's inequality, we obtain
\[\mathbb E\left(\sup_{0\le t\leq T}\|u(t)-\widetilde u(t)\|_H^2\right)+\mathbb E\left(\sup_{0\le t\leq T}\|\eta(t)-\widetilde\eta(t)\|_{M_\mu^0}^2\right)\leq C(T)\|U_0-\widetilde U_0\|_{\mathcal H_\mu^0}^2.\]
Thus, the solution map $U_0\mapsto U(\cdot;U_0)$ is continuous from $\mathcal H_\mu^0$
into $L^2\bigl(\Omega;C([0,T];\mathcal H_\mu^0)\bigr).$

\textbf{Step 5. Extension of global solution.} By repeating the above arguments, we obtain the existence of the unique
solution of \eqref{lii10} on the interval $[T,2T]$, which finally leads to the completion of
the global existence and uniqueness of solutions to \eqref{lii10} by further iterations. The proof of this theorem is completed.
\end{proof}
\section{Existence of invariant measures.}
\setcounter{equation}{0}
In this section, our goal is to establish the existence of an invariant measure
$\nu^\mu$ for the Markov semigroup $P_t^\mu$ which will be introduced below. We first recall the definition of invariant measures. Under the assumptions of the well-posedness result Theorem \ref{the1}, we introduce the Markov transition probabilities associated with
the process $U(t;U_0)$ satisfying~\eqref{lii10}.
For $t\ge 0$, $U_0\in \mathcal H^{0}_{\mu}$ and any Borel set
$A\subset \mathcal H^{0}_{\mu}$, define
\[P_t^{\mu}(U_0,A):=\mathbb P\big(U(t;U_0)\in A\big).\]
Let $\mathcal B_b(\mathcal H^{0}_{\mu})$ denote the space of bounded,
Borel-measurable functions on $\mathcal H^{0}_{\mu}$. The Markov semigroup corresponding to \eqref{lii10} is the operator,
\[P_t^{\mu}f(U_0):=\mathbb E\big[f\big(U(t;U_0)\big)\big]=\int_{\mathcal H^{0}_{\mu}}f(y)P_t^{\mu}(U_{0},\mathrm{d}y),\qquad f\in\mathcal B_b(\mathcal H^{0}_{\mu}).\]
Consequently, we are able to identify $(P_t^{\mu})^*:\mathcal{P}r(\mathcal H^{0}_{\mu})\rightarrow\mathcal{P}r(\mathcal H^{0}_{\mu})$, which implies
\[(P_t^{\mu})^*\nu(A)=\int_{\mathcal H^{0}_{\mu}} P_t^{\mu}(U_{0},A)\nu(\mathrm{d}U_{0}),\]
where $\nu\in\mathcal{P}r(\mathcal H^{0}_{\mu})$ is a probability measure. In particular, $\nu$ is called
invariant for the semigroup $\{P_t^{\mu}\}_{t\ge 0}$ if for every $f\in \mathcal B_b(\mathcal H^{0}_{\mu})$,
\[\int_{\mathcal H^{0}_{\mu}} f(U_0)\, (P_t^{\mu})^*\nu(\mathrm{d}U_0)=\int_{\mathcal H^{0}_{\mu}}f(U_0)\nu(\mathrm{d}U_0),\]
i.e.,
\[\int_{\mathcal H^{0}_{\mu}}f(U_0)(P_t^{\mu})^*\nu(\mathrm{d}U_0)=\int_{\mathcal H^{0}_{\mu}}P_t^{\mu}f(U_0)\nu(\mathrm{d}U_0).\]

Now, we present the main result concerning the existence of an invariant measure in $\mathcal{H}^0_\mu$. %In order to unify the notation used below, we introduce the operator $B$ associated with the fractional diffusion term. More precisely, define the bilinear form $b:V\times V\to\mathbb R$ by
For the higher order estimates used below,
we impose the following additional assumption on the diffusion coefficient.

$(Q_2)$ For each $k=1,2$, there exists a constant $C_{k}>0$ such that $\sigma:\mathbb{H}^{k\alpha}\to L_2(\mathcal{U},\mathbb{H}^{k\alpha})$
and
\[\left\|A_\alpha^{k/2}\sigma(u)\right\|_{\mathcal{L}_2(\mathcal{U},H)}\leq
C_{k}\left(1+\left\|A_\alpha^{k/2}u\right\|_H\right),~u\in \mathbb{H}^{k\alpha}.\]
In addition, it is necessary to enumerate the additional conditions below.
\begin{itemize}
\item[$(\Lambda_1)$]  $\lambda>\max\{C_{\sigma}^{2}-a_{2},~a_{\varphi}+C_{1}^{2}\}.$
\end{itemize}
\begin{thm}\label{the2}
Suppose assumptions of Theorem \ref{the1},  $(Q_2)$ and $(\Lambda_1)$ hold. Then,  problem \eqref{lii10} has at least one invariant measure $\nu^\mu\in$  $\mathcal{P}r(\mathcal{H}_{\mu}^{0})$. %Moreover, for any
%$\beta>0$ sufficiently small independent of $\mu$, invariant measure $\nu^\mu$ of \eqref{lii10} satisfies
%\begin{align}\label{h1}
%\sup_{\mu}\int_{\mathcal{H}^0_\mu}\exp\left(\beta\|(u,\eta)\|_{\mathcal{H}^0_\mu}^2\right)\nu^\mu(\mathrm{d}u,\mathrm{d}\eta)<\infty.
%\end{align}
\end{thm}

To establish the existence of an invariant measure to \eqref{lii10}, to apply the classical Krylov-Bogoliubov approach \cite{7} successfully, it remains to show that $P_t^\mu$ is a Feller semigroup and the family of probability measures $\{\nu_n^\mu\}_{n=1}^{\infty}$ defined by
\begin{align}\label{4.1}
\nu_n^\mu(\cdot)=\frac{1}{n}\int_0^n P_t^\mu(0,(0,0);t,(\cdot,\cdot))\mathrm{d}t,
\end{align}
is tight in $\mathcal H^{0}_{\mu}$, where  $n\in \mathbb{N}$ and $P_t^\mu(0,(0,0);t,(\cdot,\cdot))$ denotes the Markov transition probabilities associated with the process $U(t;(0,0))=(u(t;0), \eta(t,\cdot;0))$ that satisfies \eqref{lii10} starting from 0.

Note that Theorem \ref{the1} ensures that $P_t^\mu$ is Feller.  That is, $P_t^\mu f \in C_{B}(\mathcal{H}^0_\mu)$ holds for every $f\in C_{B}(\mathcal{H}^0_\mu)$, where $C_B(\mathcal H_\mu^0)$ denotes the space of all bounded continuous
real-valued functions on $\mathcal H_\mu^0$.  Hence, it is sufficient to show the probability measures family $\{\nu_n^\mu\}_{n=1}^{\infty}$ is tight.

To prove Theorem \ref{the2}, we begin with two auxiliary lemmas. Define the Lyapunov functional
\begin{equation}\label{Z6}
\Psi_0(u,\eta)=\frac{1}{2}\|u\|_H^2+\frac{1-\kappa}{2}\|\eta\|_{M^0_\mu}^2.
\end{equation}
\begin{lem}\label{Z5}
Assume hypotheses of Theorem \ref{the1}, $(Q_1)$ and $(\Lambda_1)$ hold. Let $U(t;(u_0,\eta_0)):=(u(t;u_{0}),\eta(t,\cdot;\eta_{0}))$ be the solution to system \eqref{lii10} corresponding to the initial condition $U_0=(u_0,\eta_0)\in \mathcal{H}^0_\mu$. Then
\begin{equation}\label{Y}
\mathbb{E}\Psi_0(u(t),\eta(t))\leq e^{-\hat{c}_{0}t}\Psi_0(u_0,\eta_0)+\hat{C}_{0},
\qquad t\geq0,
\end{equation}
where $\hat{c}_0,\hat{C}_{0}>0$ are independent of $t\geq0$ and $U_0$.
\end{lem}
\begin{proof}
Applying It\^o's formula to $\Psi_0(u(t),\eta(t))$ yields
\begin{equation}\label{U2}
\mathrm{d}\Psi_0(u(t),\eta(t))
=\mathcal{L}^\mu\Psi_0(u(t),\eta(t))\mathrm{d}t+\big\langle u(t),\sigma(u(t))\mathrm{d}W(t)\big\rangle_H,
\end{equation}
where $\mathcal{L}^\mu$ is the generator defined in \eqref{lili}.
Using $D_u\Psi_0=u$, $D_\eta\Psi_0=(1-\kappa)\eta$, $D_{uu}\Psi_0=I_H$, we obtain
\begin{align*}
\mathcal{L}^\mu\Psi_0(u,\eta)
=&-\kappa a_{\alpha}(u,u)+(1-\kappa)\langle \mathcal{T}_\mu\eta,\eta\rangle_{M^0_\mu}
+\langle\varphi(u),u\rangle_{\langle L^{q}(\mathcal{O}),L^{p_{0}+1}(\mathcal{O})\rangle}\nonumber\\
&+\frac12\|\sigma(u)\|_{\mathcal{L}_2(\mathcal{U},H)}^2-\lambda\|u\|_{H}^{2}.
\end{align*}
By means of \eqref{lii1}, one has
\[\langle \mathcal{T}_\mu\eta,\eta\rangle_{M^0_\mu}\leq-\frac{\delta}{2}\|\eta\|_{M^0_\mu}^2.\]
Making use of $(P_2)$, we have
\[\langle\varphi(u),u\rangle_{\langle L^{q}(\mathcal{O}),L^{p_{0}+1}(\mathcal{O})\rangle}\leq-a_2\|u\|_{L^{p_0+1}}^{p_0+1}+a_3|\mathcal O|.\]
%Moreover, by \textnormal{(Q3)} there exists $C>0$ such that
%\[\|\sigma(u)\|_{\mathcal{L}_2(U,H)}^2\leq C(1+\|u\|_H^2).\]
In particular, when $p_0=1$, by assumptions $(P_2)$ and $(Q_1)$, we have
\[\mathcal{L}^\mu \Psi_0(u,\eta)\leq -\kappa  a_{\alpha}(u,u)-\frac{(1-\kappa)\delta}{2}\|\eta\|_{M_\mu^0}^2-\bigl(\lambda+a_2-C_\sigma^2\bigr)\|u\|_H^2+C.\]
Then, by $(\Lambda_{1})$, we have
\[\mathcal{L}^\mu\Psi_0(u,\eta)\leq -\hat{c}_0\Psi_0(u,\eta)+C.\]
On the other hand, when $p_0>1$, since $\mathcal O$ is bounded, the embedding $L^{p_0+1}(\mathcal O)\hookrightarrow L^2(\mathcal O)$, together with
Young's inequality, yields
\[C_\sigma^2\|u\|_H^2\leq\frac{a_2}{2}\|u\|_{L^{p_0+1}}^{p_0+1}+C.\]
Consequently,
\[\mathcal{L}^\mu\Psi_0(u,\eta)\leq-\kappa a_{\alpha}(u,u)-\frac{(1-\kappa)\delta}{2}\|\eta\|_{M_\mu^0}^2-\lambda\|u\|_H^2-\frac{a_2}{2}\|u\|_{L^{p_0+1}}^{p_0+1}+C.\]
Recalling the definition of $\Psi_0(u,\eta)$ in \eqref{shu}, by $(\Lambda_1)$, we infer that, in both cases, there exists a constant $\hat{c}_0>0$ such that
\begin{align}\label{huai}
\mathcal{L}^\mu \Psi_0(u,\eta)\leq -\hat{c}_0\Psi_0(u,\eta)+C,
\end{align}
where
$\hat{c}_0=\begin{cases}\displaystyle\min\left\{2\bigl(\lambda+a_2-C_\sigma^2\bigr),\delta\right\},
& p_0=1,\\[2mm]
\displaystyle
\min\left\{2\lambda,\,\delta\right\},
& p_0>1.
\end{cases}
$
By taking expectation in \eqref{U2} and using \eqref{huai},  %we obtain
%\[\frac{\mathrm{d}}{\mathrm{d}t}\mathbb E\Psi_0(u(t),\eta(t))\leq-c_0\mathbb E\Psi_0(u(t),\eta(t))+C.\]
it follows from Gronwall's inequality that
\[\mathbb E\Psi_0(u(t),\eta(t))\leq e^{-\hat{c}_0 t}\Psi_0(u_0,\eta_0)+\frac{C}{\hat{c}_0}\bigl(1-e^{-\hat{c}_0 t}\bigr)\leq e^{-\hat{c}_0 t}\Psi_0(u_0,\eta_0)+\hat{C}_{0},
\qquad t\geq0.\]
This completes the proof of this lemma.
\end{proof}
Next, to derive higher order energy estimates, we can regard $A_\alpha$ as the positive self-adjoint operator on $H$ with domain $D(A_\alpha)$, as introduced in Section 2.1, rather than merely as an operator from \(V\) into \(V^*\).  We then introduce the functional below
\begin{equation}\label{huai2}
\Psi_1(u,\eta)=\frac{1}{2}\|A_{\alpha}^{1/2}u\|_{H}^2+\frac{1-\kappa}{2 }\|\eta\|_{M_\mu^\alpha}^2.
\end{equation}
\begin{lem}\label{shao}
Assume hypotheses of Theorem \ref{the1}, $(Q_2)$  and $(\Lambda_1)$ hold. Let $U_0=(u_0,\eta_0)\in\mathcal H_\mu^\alpha$,
then
\begin{equation}\label{huai5}
\mathbb E\Psi_1(u(t),\eta(t))\leq e^{-\hat{c}_{1}t}\Psi_1(u_0,\eta_0)+\hat{C}_{1},
\qquad t\geq0,
\end{equation}
where $\hat{c}_{1}, \hat{C}_{1}>0$ are independent of $U_0$ and $t$.
\end{lem}
\begin{proof}
We formally derive the uniform estimates of solutions which can be justified by a limiting process. Taking into account $D_u\Psi_{1}=A_{\alpha}u$, $D_\eta\Psi_{1}=(1-\kappa)A_{\alpha}\eta$ and $D_{uu}\Psi_{1}=A_{\alpha}I_{H}$, applying It\^o's formula to $\Psi_{1}(u(t),\eta(t))$, we have
\begin{equation}\label{longhuai6}
\mathrm{d}\Psi_{1}(u(t),\eta(t))=\mathcal L^{\mu}\Psi_{1}(u(t),\eta(t))\mathrm{d}t+\langle A_{\alpha}^{1/2}u(t),A_{\alpha}^{1/2}\sigma(u(t))\mathrm{d}W(t)\rangle_{H}.
\end{equation}
By \eqref{lili}, we obtain
\begin{align}\label{huai7}
\mathcal L^{\mu}\Psi_{1}(u,\eta)=&-\kappa\big\|A_{\alpha}u\big\|_{H}^{2}-(1-\kappa)\int_{0}^{\infty}\mu(s)\langle A_{\alpha}\eta(s),A_{\alpha}u\rangle_{H}\mathrm{d}s+(1-\kappa)\langle \mathcal{T}_{\mu}\eta,\eta\rangle_{M^{\alpha}_{\mu}}\nonumber \\
&+(1-\kappa)\int_0^\infty \mu(s)\langle A_{\alpha}\eta(s),A_{\alpha}u\rangle_H\mathrm{d}s-\lambda\langle u,A_{\alpha}u\rangle_{H}\nonumber \\
&+\big\langle \varphi(u),A_{\alpha}u\big\rangle_{H}+\frac12\big\|A_{\alpha}^{1/2}\sigma(u)\big\|^{2}_{\mathcal{L}_{2}(\mathcal{U},H)}.
\end{align}
By \eqref{2.8}, we derive
\begin{equation}\label{huai8}
\langle \mathcal{T}_\mu\eta,\eta\rangle_{M^\alpha_\mu}\leq-\frac{\delta}{2}\|\eta\|_{M^\alpha_\mu}^{2}.
\end{equation}
%Next, by $(P_0)$, $(P_3)$, the Cauchy-Schwarz  and Young inequalities, we obtain
%\begin{align*}
%\langle \varphi(u), B u\rangle_{\langle L^{q}(\mathcal{O}),L^{p_{0}+1}(\mathcal{O})\rangle}
%\leq&\|\varphi(u)\|_{H}\|Bu\|_{H}
%\leq \red{|a_\varphi| }\|u\|_{H}\|Bu\|_{H}\leq\frac{\kappa}{2}\|Bu\|_{H}^2+\frac{2a_{\varphi}^{2}}{\kappa}\|u\|_{H}^2.
%\end{align*}
%Moreover, by the fractional Poincar\'{e} inequality , we have
%\[\|u\|_{H}^2\leq \beta\|(-\Delta)^{\alpha/2}u\|_{H}^2.\]
%Consequently,
%\[|\langle \varphi(u),(-\Delta)^\alpha u\rangle_H|.\]
%%Moreover, by the fractional Poincar\'e inequality \eqref{shu1} applied to $\nabla u$,
%%there exists $\beta>0$ (independent of $\mu$) such that
%%\begin{equation}\label{huai10}
%%\big\|(-\Delta)^{\frac{\alpha+i}{2}}u\big\|_{H}^{2}\geq\beta\|u\|_{H^{i}}^{2},~~~i=1,...,m.
%%\end{equation}
For the nonlinear term, by the definition of $A_\alpha$, we have
\[\begin{aligned}
\left\langle\varphi(u),A_\alpha u\right\rangle_{H}=\frac{C(n,\alpha)}{2}(\varphi(u),u)_{\dot{H}(\mathbb{R}^{n})}=
\frac{C(n,\alpha)}{2}\int_{\mathbb R^n}\int_{\mathbb R^n}\frac{(u(x)-u(y))\bigl(\varphi(u(x))-\varphi(u(y))\bigr)}{|x-y|^{n+2\alpha}}\mathrm{d}x\mathrm{d}y.
\end{aligned}\]
By assumption $(P_3)$ and the mean value theorem, for all $r,s\in\mathbb R$,
\[\bigl(\varphi(r)-\varphi(s)\bigr)(r-s)\leq a_\varphi |r-s|^2.\]
Consequently,
\[\begin{aligned}
\left\langle\varphi(u),A_\alpha u\right\rangle_{H} \leq a_\varphi\frac{C(n,\alpha)}{2}\int_{\mathbb R^n}\int_{\mathbb R^n}\frac{|u(x)-u(y)|^2}{|x-y|^{n+2\alpha}}
\mathrm{d}x\mathrm{d}y
=a_\varphi\|A_\alpha^{1/2}u\|_H^2.
\end{aligned}
\]
Moreover, by $(Q_2)$, we find
\begin{equation}\label{huai11}
\big\|A_{\alpha}^{1/2}\sigma(u)\big\|^{2}_{\mathcal{L}_{2}(\mathcal{U},H)}
\leq 2C_{1}^{2}\left(1+\|A_{\alpha}^{1/2}u\|_{H}^{2}\right).
\end{equation}
Combining \eqref{longhuai6}-\eqref{huai11} and $(\Lambda_1)$, we infer that
\begin{align}\label{huai12}
\mathcal L^{\mu}\Psi_{1}(u,\eta)\leq&-\kappa\big\|A_{\alpha}u\big\|_{H}^{2}
-\frac{1-\kappa}{2}\delta\|\eta\|_{M^{\alpha}_{\mu}}^{2}
+(-\lambda+a_\varphi+C_{1}^{2})\|A_{\alpha}^{1/2}u\|_{H}^{2}+C_{1}^{2}\nonumber\\[1.0ex]
\leq&-\hat{c}_{1}\Psi_{1}(u,\eta)+C,
\end{align}
where  $\hat{c}_{1}=\min\left\{\delta,2(\lambda-a_\varphi-C_{1}^2)\right\}$.
Substituting \eqref{huai12} into \eqref{longhuai6} and taking expectation, we obtain
\begin{align}\label{fei1}
\frac{\mathrm{d}}{\mathrm{d}t}\mathbb{E}\Psi_{1}(u(t),\eta(t))\leq-\hat{c}_{1}\mathbb{E}\Psi_1(u(t),\eta(t))+C.
\end{align}
%By Lemma \ref{Z5}, we have
%\[\frac{\mathrm{d}}{\mathrm{d}t}\E\Psi_1(u(t),\eta(t))\leq-c\E\Psi_1(u(t),\eta(t))+Ce^{-c_0 t}\Psi_0(u_0,\eta_0)+C(C_0+1).\]
Multiplying  both sides of the above inequality by $e^{\hat{c}_{1}t}$ and integrating it over $[0,t]$, we prove
%\begin{equation}\label{liang4}
%\E\Psi_1(u(t),\eta(t))\leq e^{-ct}\Psi_1(u_0,\eta_0)+C e^{-ct}\Psi_0(u_0,\eta_0)\int_0^t e^{(c-c_0)s}\mathrm{d}s+C.
%\end{equation}
%By the fractional Poincar\'{e} inequality (cf. \eqref{shu1}) and $M_{\mu}^{\alpha}\hookrightarrow M_{\mu}^{0}$, we obtain
%\begin{align*}
%\Psi_0(u_0,\eta_0)\leq C\Psi_1(u_0,\eta_0).
%\end{align*}
%\red{In addition, it follows from Lemma 3.2 in \cite{er1} that
%$$e^{-ct}\int_{0}^{t}e^{-(c_{0}-c)r}\mathrm{d}r\leq Ce^{-\tilde{c}t},~~~t\geq0,$$}
%where $C=C(c_{0},c)$, $\tilde{c}=\tilde{c}(c_{0},c)$. Substituting the above inequality  into \eqref{liang4} yields
%\[\mathbb{E}\Psi_1(u(t),\eta(t))\leq e^{-c_1 t}\Psi_1(u_0,\eta_0)+C_1,\qquad t\geq0,\]
%for some constants $c_1>0$ and $C_1>0$ independent of $t$, $U_0$ and $\mu$.
%In view of Lemma \ref{Z5} and the Sobolev embedding, we obtain
%\[\mathbb{E}\|u(r)\|_H^2\leq2\mathbb{E}\Psi_0(U(r))\leq2 e^{-c_0 r}\Psi_0(U_0)+2C_{0}\leq 2\alpha_1^{-1}e^{-c_0 r}\Psi_1(U_0)+2C_{0},\]
%where $\alpha_1>0$ denotes the constant appearing in the Sobolev embedding.
%Consequently,
%\[\int_0^te^{-c(t-r)}\mathbb{E}\|u(r)\|_H^2\mathrm{d}r\leq 2\alpha_1^{-1}e^{-ct}\int_0^te^{-(c_0 - c) r}\mathrm{d}r\Psi_1(U_0)+C.\]
%Applying inequality $e^{-ct}\int_{0}^{t}e^{-(c'-c)r}\mathrm{d}r\leq Ce^{-\tilde{c}t},\qquad t\geq0,$
%for some positive constants $C,\tilde{c}>0$ depending only on $c,c'>0$ and independent of $t$, it follows that
%\begin{align}\label{mi1}
%e^{-ct}\int_0^t e^{-(c_0 - c)r}\mathrm{d}r\leq C e^{-ct}.
%\end{align}
 \eqref{huai5}.
\end{proof}
\begin{lem}\label{h}
Assume hypotheses of Theorem \ref{the1} and $(\Lambda_{1})$ hold. Let $U(t;(u_{0},\eta_{0}))=(u(t;u_{0}),\eta(t,\cdot;\eta_{0}))$ denote the solution of system \eqref{lii10} corresponding to the initial condition
$U_0=(u_{0},\eta_{0})\in \mathcal{H}^0_\mu$. Then, for all $0<\theta\leq\min\{\delta,2(\lambda+a_{2}-C_{\sigma}^{2})\}$, where $\delta$ is the constant  in \eqref{lii2}, there exists a positive constant $c$ independent of $t$ and $U_{0}$, such that
\begin{equation}\label{h5}
\mathbb{E}\int_0^t e^{\theta r}a_{\alpha}(u(r),u(r))\mathrm{d}r\leq c\|(u_{0},\eta_{0})\|_{\mathcal{H}^0_\mu}^2 + c e^{\theta t}.
\end{equation}
\end{lem}
\begin{proof}
Applying It\^{o}'s formula to $e^{\theta t}\Psi_0(u(t),\eta(t))$, we have
\begin{align*}
\mathrm{d}(e^{\theta t}\Psi_0(u(t),\eta(t)))=&e^{\theta t}\Bigl(\theta\Psi_0(u(t),\eta(t))-\kappa a_{\alpha}(u(t),u(t))
+(1-\kappa)\langle \mathcal{T}_{\mu}\eta(t),\eta(t)\rangle_{M_\mu^0}\\[0.6ex]
&+\langle\varphi(u(t)),u(t)\rangle_{\langle L^{q}(\mathcal{O}),L^{p_{0}+1}(\mathcal{O})\rangle}+\tfrac12\|\sigma(u(t))\|_{\mathcal{L}_2(\mathcal{U},H)}^2
-\lambda\|u(t)\|_{H}^{2}\Bigr)\mathrm{d}t\\[0.6ex]
&+e^{\theta t}\langle u(t),\sigma(u(t))\mathrm{d}W(t)\rangle_H.
\end{align*}
Using  \eqref{lii1} and the assumption $\theta\leq\delta$, we obtain
\begin{equation}\label{4271}
\tfrac{\theta(1-\kappa)}{2}\|\eta(t)\|_{M^0_\mu}^2+(1-\kappa)\langle \mathcal{T}_\mu\eta(t),\eta(t)\rangle_{M^0_\mu}\leq0.
\end{equation}
Making use of   $(P_2)$ and $(Q_1)$, along with the H\"{o}lder inequality, we obtain
\[\begin{aligned}\label{4.16}
&\frac{\theta}{2}\|u(t)\|_H^2+\langle \phi(u(t)),u(t)\rangle_{\langle L^q(\mathcal O),L^{p_0+1}(\mathcal O)\rangle}+\frac12\|\sigma(u(t))\|_{\mathcal{L}_2(U,H)}^2\\
\leq&C_\sigma^2+\left(C_\sigma^2+\frac{\theta}{2}-\lambda\right)\|u(t)\|_H^2-a_2\|u(t)\|_{L^{p_0+1}}^{p_0+1}+a_3|\mathcal O|.
\end{aligned}
\]
We distinguish the following two cases. If $p_0=1$, then $L^{p_0+1}(\mathcal O)=L^2(\mathcal O)=H$, and hence
\[\begin{aligned}\label{4.17}
&\frac{\theta}{2}\|u(t)\|_H^2+\langle \phi(u(t)),u(t)\rangle_H+\frac12\|\sigma(u(t))\|_{\mathcal{L}_2(U,H)}^2\\
\leq&C_\sigma^2+a_3|\mathcal O|-\left(\lambda+a_2-C_\sigma^2-\frac{\theta}{2}\right)\|u(t)\|_H^2.
\end{aligned}
\]
By $(\Lambda_{1})$, we obtain
\begin{align}\label{4.18}
\frac{\theta}{2}\|u(t)\|_H^2+\langle \phi(u(t)),u(t)\rangle_H+\frac12\|\sigma(u(t))\|_{\mathcal{L}_2(U,H)}^2\leq C.
\end{align}
On the other hand, if $p_0>1$, since $\mathcal O$ is bounded, the embedding $L^{p_0+1}(\mathcal O)\hookrightarrow L^2(\mathcal O)$
and Young's inequality imply that, for every fixed $\theta>0$,
\[\left(C_\sigma^2+\frac{\theta}{2}\right)\|u(t)\|_H^2\leq\frac{a_2}{2}\|u(t)\|_{L^{p_0+1}}^{p_0+1}+C.\]
Consequently,
\begin{align}\label{4.19}
&\frac{\theta}{2}\|u(t)\|_H^2+\langle \phi(u(t)),u(t)\rangle_{\langle L^q(\mathcal O),L^{p_0+1}(\mathcal O)\rangle}+\frac12\|\sigma(u(t))\|_{L_2(\mathcal{U},H)}^2\nonumber\\
\leq&-\lambda\|u(t)\|_H^2-\frac{a_2}{2}\|u(t)\|_{L^{p_0+1}}^{p_0+1}+C\leq C.
\end{align}
Thus, in both cases, for a sufficiently small $\theta>0$ with
$\theta\leq\delta$, we have
\begin{align}\label{4.20}
\frac{\theta}{2}\|u(t)\|_H^2+\langle \phi(u(t)),u(t)\rangle_{\langle L^q(\mathcal O),L^{p_0+1}(\mathcal O)\rangle}
+\frac12\|\sigma(u(t))\|_{L_2(\mathcal{U},H)}^2\leq C.
\end{align}
It follows from \eqref{4271}-\eqref{4.20} that there exists a positive constant $c=c(\theta,p_{0},a_{2},C_{\sigma},\mathcal{O})$ which is independent of $t$, such that
\begin{equation*}
\begin{split}
&~\quad e^{\theta t}\mathbb{E}\Psi_0(u(t),\eta(t))+\kappa\mathbb{E}\int_0^t e^{\theta r}a_{\alpha}(u(r),u(r))\mathrm{d}r\\
&\leq c\|(u_{0},\eta_{0})\|^2_{\mathcal{H}^0_\mu}+c\int_0^te^{\theta r}\mathrm{d}r
\leq c\|(u_{0},\eta_{0})\|_{\mathcal{H}^0_\mu}^2+ce^{\theta t}.
\end{split}
\end{equation*}
Thus, the proof of this lemma is complete.
\end{proof}
To derive the estimates to the solution of system \eqref{lii10} in $E_{\mu}^{\alpha}$, based on the definition of the norm of $E_{\mu}^{\alpha}$, it remains  to prove the solution of \eqref{lii10} is uniformly bounded in  $\mathcal{T}_{\mu}$ and $\mathbb{T}^{\mu}$  with respect to $t$.
\begin{lem}\label{hx2}
Assume hypotheses of Theorem \ref{the1} hold. Let $U(t;(0,0))=(u(t;0),\eta(t,\cdot;0))$ be the solution of system \eqref{lii10} with zero initial data.
Then, there exists a constant $C$ such that the following estimates hold for all $t \geq0$,
\[\mathbb{E}\|\mathcal{T}_{\mu}\eta(t)\|^2_{M^0_{\mu}}\leq C\qquad\text{and}\qquad\mathbb{E}\sup_{r \geq1} r\mathbb{T}^{\mu}_{\eta(t)}(r)\leq C.\]
\end{lem}
\begin{proof}
Set $\theta_0:=\min\left\{\frac{\delta}{2},\lambda+a_2-C_\sigma^2\right\}>0.$ By $(\Lambda_1)$, $\theta_0$ is well defined, and moreover
$0<\theta_0\leq\min\left\{\delta,2(\lambda+a_2-C_\sigma^2)\right\}.$
Since $\eta_0(s)=0$, by \eqref{x}, we find
\begin{align}\label{zuihou}
\eta(t,s)=\chi_{(0,t]}(s)\int_{0}^{s}u(t-r)\mathrm{d}r+\chi_{(t,\infty)}(s)\int_{0}^{t}u(t-r)\mathrm{d}r,
\end{align}
which implies
\[\mathcal{T}_{\mu}\eta(t,s)=-u(t-s)\chi_{(0,t]}(s).\]
Thus, it is obvious that
\[\mathbb{E}\|\mathcal{T}_{\mu}\eta(t)\|^2_{M^0_{\mu}}=\mathbb{E}\int_{0}^{t} \mu(s)a_{\alpha}(u(t-s),u(t-s))\mathrm{d}s
=\mathbb{E}\int_{0}^{t} \mu(t-s)a_{\alpha}(u(s),u(s))\mathrm{d}s.\]
Using  \eqref{lii2} and Lemma \ref{h}, we obtain
\begin{align*}
\mathbb{E}\|\mathcal{T}_{\mu}\eta(t)\|^2_{M^0_{\mu}}\leq&
\mu(0)e^{-\delta t}\mathbb{E}\int_{0}^{t}e^{\delta s}a_{\alpha}(u(s),u(s))\mathrm{d}s\nonumber\\
\leq&\mu(0)e^{(\delta-\theta-\delta) t}\mathbb{E}\int_{0}^{t}e^{\theta s}a_{\alpha}(u(s),u(s))\mathrm{d}s\nonumber\\
\leq&\mu(0)e^{-\theta t}\left(\|U_{0}\|_{\mathcal{H}^0_{\mu}}^{2}+ce^{\theta t}\right)\leq C.
\end{align*}
Similarly, for $\sup_{r\ge1} r\mathbb{T}_{\eta(t)}^\mu(r)$, taking into account the representation formula of $\eta$ (cf. \eqref{zuihou}) and using the Cauchy-Schwarz inequality, we infer
\begin{align} \label{o1}
a_{\alpha}(\eta(t,s),\eta(t,s))&=
a_{\alpha}(\int_0^{s\wedge t}u(t-r_1)\mathrm{d}r_1,\int_0^{s\wedge t}u(t-r_2)\mathrm{d}r_2) \nonumber\\
%&=\int_0^{s\wedge t}\int_0^{s\wedge t}
%\langle B^{1/2}u(t-r_1),B^{1/2}u(t-r_2)\rangle_H\mathrm{d}r_1\mathrm{d}r_2 \nonumber\\
&\leq(s\wedge t)\int_0^{s\wedge t}a_{\alpha}(u(t-r'),u(t-r'))\mathrm{d}r' \nonumber\\
&\leq s\int_0^{s\wedge t}a_{\alpha}(u(t-r'),u(t-r'))\mathrm{d}r'.
\end{align}
For $s\in[0,t]$, let us define the auxiliary function
\begin{equation}\label{o2}
I(t\wedge s):=e^{-\theta_{0}s}\int_0^s a_{\alpha}(u(t-r'),u(t-r'))\mathrm{d}r'.
\end{equation}
For future use, we recall the elementary estimate that for any nonnegative locally integrable function $f:\mathbb{R}\to[0,\infty)$ with  $0<s\leq t$ and  $\delta>0$,
\begin{align}\label{min}
e^{-\theta_{0}s}\int_0^s f(t-r)\mathrm{d}r\leq\int_0^t e^{-\theta_{0}r}f(t-r)\mathrm{d}r.
\end{align}
Combining \eqref{lii2}, \eqref{2.9} with \eqref{o1}, for all $r\geq1$, we have
\begin{align}\label{o3}
r\mathbb{T}_{\eta(t)}^\mu(r)&=r\int_{(0,1/r)\cup(r,\infty)}\mu(s)a_{\alpha}(\eta(t,s),\eta(t,s))\mathrm{d}s \nonumber\\
&\leq r\mu(0)\int_{(0,1/r)\cup(r,\infty)}se^{-\delta s}\int_0^{s\wedge t}a_{\alpha}(u(t-r'),u(t-r'))\mathrm{d}r'\mathrm{d}s\nonumber\\
&\leq r\mu(0)\int_{(0,1/r)\cup(r,\infty)}s e^{-\theta_{0}s}
\left(e^{-\theta_{0}(s\wedge t)}\int_0^{s\wedge t}a_{\alpha}(u(t-r'),u(t-r'))\mathrm{d}r'\right)\mathrm{d}s\nonumber\\
&\leq r\mu(0)\int_{(0,1/r)\cup(r,\infty)} s e^{-\theta_{0}s}I(s\wedge t)\mathrm{d}s.
\end{align}
By definition \eqref{o2} and \eqref{min}, we find
\begin{equation}\label{O}
I(t\wedge s)\leq\int_0^te^{-\theta_{0}r'}a_{\alpha}(u(t-r'),u(t-r'))\mathrm{d}r'.
\end{equation}
Therefore,
\begin{align}\label{O1}
r\mathbb{T}_{\eta(t)}^\mu(r)&\leq r\mu(0)\left(\int_{(0,1/r)\cup(r,\infty)}se^{-\theta_{0}s}\mathrm{d}s\right)
\left(\int_0^te^{-\theta_{0}r'}a_{\alpha}(u(t-r'),u(t-r'))\mathrm{d}r'\right).
\end{align}
For $r\geq1$, by integration by parts, we obtain
\[r\int_{(0,1/r)\cup(r,\infty)}se^{-\theta_{0}s}\mathrm{d}s
<C(\delta)<\infty.\]
Taking expectation and supremum in \eqref{O1}, using the a priori energy estimate \eqref{h5}, we conclude that
\begin{align*}
\mathbb{E}\sup_{r\geq1}r\mathbb{T}_{\eta(t)}^\mu(r)&\leq\mathbb{E}\left(\sup_{r\geq1} r\mu(0)
\int_{(0,\frac{1}{r})\cup(r,\infty)}se^{-\theta_{0}s}\mathrm{d}s
\int_0^t e^{-\theta_{0}r'}a_{\alpha}(u(t-r'),u(t-r'))\mathrm{d}r'\right)\\
&\leq C(\delta)\mu(0)\mathbb{E}\int_0^t e^{-\theta_{0}r'}a_{\alpha}(u(t-r'),u(t-r'))\mathrm{d}r'\\
&=C(\delta)\mu(0)e^{-\theta_{0}t}
\mathbb{E}\int_0^te^{\theta_{0}s}a_{\alpha}(u(s),u(s))\mathrm{d}s\leq C.
\end{align*}
 This completes the proof.
\end{proof}
We next prove Theorem \ref{the2} to show  the existence of an invariant measure
$\nu^{\mu}$ to \eqref{lii10}.
\begin{proof}
We prove existence of an invariant measure  to problem \eqref{lii10} by  Krylov-Bogoliubov method.
 It follows from Theorem \ref{the1}  that the semigroup $\{P_t^\mu\}_{t\geq0}$ is Feller, by \eqref{4.1} it suffices to show that $\{\nu_n^\mu\}_{n=1}^{\infty}$ is tight in $\mathcal{H}_\mu^0$. However, it is worth noticing
$M_{\mu}^{\alpha}\hookrightarrow M_{\mu}^{0}$ is only continuous, so is the embedding $\mathcal{H}^\alpha_\mu \hookrightarrow\mathcal{H}^0_\mu$. In contrast, the embedding $\mathcal{Z}^\alpha_\mu =\mathbb{H}^\alpha \times E^\alpha_\mu \hookrightarrow\mathcal{H}^0_\mu=H\times M_\mu^0$ is compact.

By \eqref{huai5}, together with Lemma \ref{hx2}, we deduce that there exists $\tilde{C}>0$ such that for all $t\geq0$,
\begin{align*}
\int_0^t\mathbb{E}\left(\|A_{\alpha}^{1/2}u(r)\|_H^2+\|\eta(r)\|_{E_\mu^\alpha}^2\right)\mathrm{d}r\leq \tilde{C}t,
\end{align*}
and hence for all $k\in \mathbb{N}$,
\begin{align*}
\frac{1}{k}\int_0^k\mathbb{E}\left(\|A_{\alpha}^{1/2}u(r)\|_H^2+\|\eta(r)\|_{E_\mu^\alpha}^2\right)\mathrm{d}r \leq \tilde{C}.
\end{align*}
Furthermore, for any $\delta>0$, there exists a constant $M=M(\delta)>0$ such that $\frac{\tilde{C}}{M}<\delta$. Define
$$Y^{\delta}=\{(u,\eta)\in \mathcal{Z}_\mu^\alpha:\|(u,\eta)\|^{2}_{\mathcal{Z}_\mu^\alpha}\leq M(\delta)\},$$
then the closure of $Y^{\delta}$ is compact in $\mathcal{H}_\mu^0$. By Chebyshev's inequality, we obtain
\begin{equation*}
\nu_n^\mu(Y^{\delta})\geq1-\frac{\frac{1}{n}\displaystyle\int_{0}^{n}\mathbb{E}\left(\|A_{\alpha}^{1/2}u(r)\|_H^2+\|\eta(r)\|_{E_\mu^\alpha}^2\right)\mathrm{d}r}{M}>1-\delta.
\end{equation*}
Hence $\{\nu_n^\mu\}_{n=1}^{\infty}$ is tight in $\mathcal{H}_\mu^0$.
By the Krylov-Bogoliubov theorem,  we know $\nu^\mu\in\mathcal{P}r(\mathcal{H}_\mu^0)$ is invariant
for \eqref{lii10}. This completes the proof of Theorem \ref{the2}.
\end{proof}
\begin{cor}\label{qi}
Suppose assumptions of Theorem \ref{the2} hold. Let $\nu^\mu$ be an invariant measure on $\mathcal{P}r(\mathcal{H}_{\mu}^{0})$ obtained as a weak limit of the Krylov-Bogoliubov averages defined in \eqref{4.1}. Then, there exists a constant $\hat{C}_{0}>0$ ($\hat{C}_{0}$ is the same as in Lemma \ref{Z5}) such that
\begin{equation}\label{qi1}
\int_{\mathcal{H}_\mu^0}\Psi_0(u,\eta)\nu^\mu(\mathrm{d}u,\mathrm{d}\eta)\leq \hat{C}_{0}.
\end{equation}
Moreover, there exists a constant $C>0$ such that
\begin{equation}\label{qi2}
\int_{\mathcal{H}_\mu^0}\|(u,\eta)\|_{\mathcal{H}_\mu^0}^2\nu^\mu(\mathrm{d}u,\mathrm{d}\eta)\leq C.
\end{equation}
\end{cor}
\begin{proof}
Let $\nu_n^\mu(\cdot)$
be a family of measures defined by \eqref{4.1}. By Lemma \ref{Z5}, we have
\begin{align*}
\int_{\mathcal{H}_\mu^0}\Psi_0(u,\eta)\nu_n^\mu(\mathrm{d}u,\mathrm{d}\eta)
&=\frac1n\int_0^n \int_{\mathcal{H}_\mu^0}\Psi_0(u,\eta)P_t^\mu(0,(0,0);t,(\mathrm{d}u,\mathrm{d}\eta))\mathrm{d}t\\
&=\frac1n\int_0^n \mathbb{E}\Psi_0((u(t;0),\eta(t,\cdot;0)))\mathrm{d}t\\
&\leq\frac1n\int_0^n \bigl(e^{-\hat{c}_0t}\Psi_0(0)+\hat{C}_0\bigr)\mathrm{d}t
=\hat{C}_0.
\end{align*}
Therefore,
\[\sup_{n\in\mathbb{N}}\int_{\mathcal{H}_\mu^0}\Psi_0(u,\eta)\nu_n^\mu(\mathrm{d}u,\mathrm{d}\eta)\leq \hat{C}_0.\]

By the proof of Theorem  \ref{the2}, the family
$\{\nu_n^\mu\}_{n=1}^{\infty}$ is tight in $\mathcal H_\mu^0$.
Then, there exist a probability measure
$\bar\nu^\mu \in \mathcal{P}r(\mathcal H_\mu^0)$ and a subsequence,
still denoted by $\{\nu_n^\mu\}_{n=1}^{\infty}$, such that
\[\nu_n^\mu \xrightarrow{\mathrm{weakly}} \bar\nu^\mu,\qquad n\to\infty .\]
Moreover, by the Krylov--Bogoliubov argument in the proof of Theorem \ref{the2},
$\bar\nu^\mu$ is an invariant measure.
For simplicity, we still denote $\bar\nu^\mu$ by $\nu^\mu$.

Furthermore, we have
\[\int_{\mathcal{H}_\mu^0}\Psi_0(u,\eta)\nu^\mu(\mathrm{d}u,\mathrm{d}\eta)\leq\liminf_{n\to\infty}\int_{\mathcal{H}_\mu^0}\Psi_0(u,\eta)\nu_{n}^\mu(\mathrm{d}u,\mathrm{d}\eta)\leq \hat{C}_0.\]
This proves \eqref{qi1}.  Finally, by the definition of $\Psi_0$, we know
\[\|(u,\eta)\|_{\mathcal{H}_\mu^0}^2\leq\frac{2}{1-\kappa}\Psi_0(u,\eta), \qquad (u,\eta)\in \mathcal{H}_\mu^0.\]
Combining  with \eqref{qi1}, we obtain
\[\int_{\mathcal{H}_\mu^0}\|(u,\eta)\|_{\mathcal{H}_\mu^0}^2\nu^\mu(\mathrm{d}u,\mathrm{d}\eta)
\leq \frac{2}{1-\kappa}\int_{\mathcal{H}_\mu^0}\Psi_0(u,\eta)\nu^\mu(\mathrm{d}u,\mathrm{d}\eta)\leq C.\]
The proof of this corollary is finished.
\end{proof}
\section{Regularity of the invariant measure $\nu^{\mu}$.}
\setcounter{equation}{0}
In this section, we show that the support of invariant measure $\nu^{\mu}$ is contained in $\mathbb{H}^{2\alpha}\times M_{\mu}^{2\alpha}$.
We now state the main result of this section. %\red{First, we impose the following assumption:
To establish the higher-order regularity estimates, we further assume that the nonlinear term satisfies the following condition.
\begin{enumerate}
\item[($P_4$)] For every $\varepsilon>0$, there exists a constant $C_{\varepsilon,2}>0$ such that
$$\bigl(A_\alpha^{1/2}\varphi(u),
A_\alpha^{3/2}u\bigr)_H
\leq \varepsilon\left\|A_\alpha^{3/2}u\right\|_H^2+C_{\varepsilon,2}\left(1+\left\|A_\alpha u\right\|_H^2\right).$$
\item[$(\Lambda_2)$]  $\lambda>\{C_{\frac{\kappa}{2},2}+C_{2}^{2}\}.$
\end{enumerate}
%\begin{itemize}
%\item[$(P_4)$]
%$\varphi : H^{(m-1)\alpha} \to H^{(m-1)\alpha},~~~~m\geq 2.$
%\end{itemize}}
\begin{theorem}\label{chen}
Assume the conditions in Theorem \ref{the2}, $(P_4)$  and $(\Lambda_2)$ hold. Then the invariant measure $\nu^\mu$  obtained in Theorem \ref{the2} satisfies
\[
\nu^\mu(\mathbb{H}^{2\alpha}\times M^{2\alpha}_\mu)=1.
\]
\end{theorem}
Since the history variable $\eta$ is governed by a hyperbolic-type equation, no additional regularity can be obtained from its evolution alone. When proving the existence of an invariant measure, we set the initial data for the solution to be $U_{0}=(u_{0},\eta_{0})\in\mathcal{H}_{\mu}^{0}$. Hence, the regularity of an invariant measure cannot be determined by direct computation.  To overcome this difficulty,
we introduce a suitable control argument, which allows us to obtain the required
regularizing properties and asymptotically embed the dynamics into a more
regular subspace \cite{er1}. First, let us introduce the following controlled system,
\begin{equation}\label{chen1}
\left\{
\begin{aligned}
&\mathrm{d}\hat{u}(t)=
-\kappa A_{\alpha}\hat{u}(t)\mathrm{d}t
-(1-\kappa)\int_{0}^{\infty}\mu(s)A_{\alpha}\hat{\eta}(t,s)\mathrm{d}s\mathrm{d}t
+\varphi(\hat{u}(t))\mathrm{d}t\\[0.8ex]
&+\sigma(\hat{u}(t))\mathrm{d}W(t)-\kappa \lambda_{\hat{n}}\Pi_{\hat{n}}(\hat{u}(t)-u(t))\mathrm{d}t-\lambda\hat{u}(t)\mathrm{d}t,\qquad (t,x)\in \mathbb{R}^{+}\times \mathcal{O},\\[0.8ex]
&\frac{\mathrm{d}}{\mathrm{d}t}\hat{\eta}(t)=\mathcal{T}_\mu\hat{\eta}(t)+\hat{u}(t),\qquad (t,x)\in \mathbb{R}^{+}\times \mathcal{O},\\[0.8ex]
&\left(\hat{u}(0),\hat{\eta}(0)\right)=\left(0,0\right)\in\mathcal{H}_{\mu}^{0},\qquad x\in\mathcal{O},
\end{aligned}
\right.
\end{equation}
where for each $N\in\mathbb N$, define the spectral projection $\Pi_N:H\to\operatorname{span}\{f_1,\ldots,f_N\}$ by
\[\Pi_Nu=\sum_{j=1}^N\langle u,f_j\rangle_Hf_j,\qquad u\in H,\]
and $\{f_j\}_{j=1}^{\infty}$ is the orthonormal basis of eigenfunctions of $A_\alpha$ introduced in Section 2.1, satisfying $A_\alpha f_j=\lambda_jf_j$.
 Choose $\hat n$ sufficiently large so that $\lambda_{\hat n}$ satisfies,
\begin{equation}\label{shui}
\lambda+\kappa\lambda_{\hat n}>a_{\varphi}+\frac{1}{2}L_{\sigma}^{2}.
\end{equation}
Next, we present several lemmas that will be used in the proof of the main result.

Initially, we focus on deriving moment estimates
for the solution in higher-order regularity under the additional
assumption $(Q_2)$, $(P_4)$ and $(\Lambda_2)$. Define
\begin{equation}\label{huai3}
\Psi_2(u,\eta)=\frac{1}{2}\|A_{\alpha}u\|_{H}^2+\frac {1-\kappa}{2}\|\eta\|_{M_\mu^{2\alpha}}^2.
\end{equation}
\begin{lemma}\label{xiao}
Let $U_{0}=(u_{0},\eta_{0})\in \mathcal{H}_{\mu}^{2\alpha}$. Suppose  assumptions of Theorem \ref{the1}
hold true.  In addition, let  $(Q_2)$, $(P_4)$  and $(\Lambda_2)$ hold.
Then, there exist positive constants
$\hat{c}_{2}$ and $\hat{C}_{2}$, independent of $U_0$ and $t$, such that
\begin{equation}\label{xiao1}
\mathbb E\Psi_2(u(t),\eta(t))\leq e^{-\hat{c}_{2}t}\Psi_2(u_{0},\eta_{0})+\hat{C}_{2}, ~~~t\geq0.
\end{equation}
\end{lemma}
\begin{proof}
We first derive the higher order energy estimates formally, which can be justified by Galerkin approximations and a limiting argument.
For $(u,\eta)\in \mathcal{H}_\mu^{2\alpha}$, by \eqref{huai3} and direct computations, we derive
\[D_u\Psi_2=A_{\alpha}^{2}u,\qquad D_\eta\Psi_2=(1-\kappa)A_{\alpha}^{2}\eta,\qquad D_{uu}\Psi_2=A_{\alpha}^{2}I_{H}.\]
Applying It\^o's formula to $\Psi_2(u(t),\eta(t))$, we obtain
\begin{equation}\label{liang11}
\mathrm{d}\Psi_2(u(t),\eta(t))=\mathcal{L}^\mu\Psi_2(u(t),\eta(t))\mathrm{d}t+\bigl\langle A_{\alpha}u(t),A_{\alpha}\sigma(u(t))\mathrm{d}W(t)\bigr\rangle_H.
\end{equation}
Hence, by \eqref{lili}, we have
\begin{align}\label{yan}
\mathcal{L}^\mu\Psi_2(u,\eta)
=&-\kappa\|A_{\alpha}^{3/2}u\|_H^2-(1-\kappa)\int_0^\infty \mu(s)
\bigl\langle A_{\alpha}^{3/2}\eta(s),A_{\alpha}^{3/2}u\bigr\rangle_H\mathrm{d}s \\[0.8ex]
&+(1-\kappa)\langle \mathcal{T}_\mu\eta,\eta\rangle_{M_\mu^{2\alpha}}
+(1-\kappa)\int_0^\infty \mu(s)\bigl\langle A_{\alpha}^{3/2}\eta(s),A_{\alpha}^{3/2}u\bigr\rangle_H\mathrm{d}s \notag\\[0.8ex]
&+\bigl(A_\alpha^{1/2}\varphi(u),
A_\alpha^{3/2}u\bigr)_H+\frac12\|A_{\alpha}\sigma(u)\|_{\mathcal{L}_2(\mathcal{U},H)}^2
-\lambda\|A_{\alpha}u\|_H^2 .\notag
\end{align}
Next, we deal with each term in \eqref{yan}. For the third term, it follows from  \eqref{2.8} that,
\begin{equation}\label{yan1}
\langle \mathcal{T}_\mu\eta,\eta\rangle_{M_\mu^{2\alpha}}\leq-\frac{\delta}{2}\|\eta\|_{M_\mu^{2\alpha}}^2.
\end{equation}
For the nonlinear term, by assumption $(P_4)$ with
$\varepsilon=\frac{\kappa}{2}$, we have
\begin{equation}
\begin{aligned}\label{yan2}
\bigl(A_\alpha^{1/2}\varphi(u),
A_\alpha^{3/2}u\bigr)_H\leq
\frac{\kappa}{2}\left\|A_\alpha^{3/2}u\right\|_H^2+ C_{\frac{\kappa}{2},2}\left(1+\left\|A_\alpha u\right\|_H^2\right).
\end{aligned}
\end{equation}
For the sixth term, by  assumption $(Q_2)$, we have
\begin{equation}\label{yan3}
\frac12\|A_{\alpha}\sigma(u)\|_{\mathcal{L}_2(\mathcal{U},H)}^2
\leq C_{2}^2+C_{2}^2\|A_{\alpha}u\|_H^2.
\end{equation}
Substituting \eqref{yan1}--\eqref{yan3} into \eqref{yan}, we infer that
\begin{align}\label{yan4}
\mathcal{L}_\mu\Psi_2(u,\eta)
\leq&-\frac{\kappa}{2}\|A_{\alpha}^{3/2}u\|_H^2-\frac{(1-\kappa)\delta}{2}\|\eta\|_{M_\mu^{2\alpha}}^2
+\bigl(C_{2}^2+C_{\frac{\kappa}{2},2}-\lambda\bigr)\|A_{\alpha}u\|_H^2+C_{2}^2+C_{\frac{\kappa}{2},2} .
\end{align}
By $(\Lambda_2)$, we conclude that there exist
constants $\hat{c}_{2}>0$ and $C>0$, independent of $U_0$ and $t$, such that
\begin{equation}\label{yan5}
\mathcal{L}_\mu\Psi_2(u,\eta)
\leq-\hat{c}_{2}\Psi_2(u,\eta)+C .
\end{equation}
Combining \eqref{liang11} and \eqref{yan5}, and then taking expectation, we obtain
%Therefore, we have
%\begin{equation}\label{yan7}
%\mathbb E\Psi_m(U(t))\leq e^{-ct}\Psi_m(U_0)+C\int_0^t e^{-c(t-r)}\mathbb E\Psi_{m-1}(U(r))\mathrm{d}r+C.
%\end{equation}
%By means of  \eqref{liang10} and    \cite[Lemma 3.2]{er1}, we find
%\begin{align}\label{yan8}
%\int_0^t e^{-c(t-r)}\mathbb E\Psi_{m-1}(U(r))\mathrm{d}r
%\leq&C_{m-1}\int_0^t e^{-c(t-r)}\Bigl(e^{-c_{m-1}r}\Psi_{m-1}(U_0)+1\Bigr)\mathrm{d}r\nonumber\\
%\leq&C\int_0^t e^{-c(t-r)}\Bigl(e^{-c_{m-1}r}\Psi_{m}(U_0)+1\Bigr)\mathrm{d}r\nonumber\\
%\leq& Ce^{-ct}\Psi_{m}(U_0)+C.
%\end{align}
\[\mathbb E\Psi_2(u(t),\eta(t))
\leq e^{-\hat{c}_2 t}\Psi_2(u_0,\eta_0)+\hat{C}_2,
\qquad t\geq0,\]
for some positive constants $\hat{c}_2$ and $\hat{C}_2$ independent of $U_0$ and $t$. Thus, the proof of this lemma is complete.
\end{proof}
\begin{lemma}\label{chen4}
Assume  conditions in Theorem \ref{the1} hold. Moreover, suppose that $(\Lambda_1)$, $(Q_2)$, $(P_4)$  and $(\Lambda_2)$ hold. Then, the following facts hold:
\begin{enumerate}
\item There exist constants $c,C>0$ such that, for every $U_0=(u_{0},\eta_{0})\in\mathcal{H}^0_\mu$ and $t\geq0$,
\begin{equation}\label{chen5}
\mathbb{E}\|(u(t)-\hat{u}(t),\eta(t)-\hat{\eta}(t))\|^2_{\mathcal{H}^0_\mu}
\leq Ce^{-ct}\|(u_{0},\eta_{0})\|^2_{\mathcal{H}^0_\mu}.
\end{equation}
\item There exists a positive constant $C$, independent of $U_0$ and $t$, such that,
\begin{equation}\label{chen6}
\mathbb{E}\|(\hat{u}(t),\hat{\eta}(t))\|^2_{\mathcal{H}^{2\alpha}_\mu}
 \leq C\big(\|(u_0,\eta_0)\|^{2}_{\mathcal{H}^0_\mu}+1\big).
\end{equation}
\end{enumerate}
\end{lemma}
\begin{proof}
Let $z=\hat{u}-u~\text{and}~\zeta=\hat{\eta}-\eta$ satisfy
\begin{equation}\label{zhong}
\left\{
\begin{aligned}
&\mathrm{d}z(t)=\left(-\kappa A_{\alpha}z(t)-(1-\kappa)\int_0^{\infty} \mu(s)A_{\alpha}\zeta(t,s)\mathrm{d}s+\varphi(\hat{u}(t))-\varphi(u(t))\right.\\[1.0ex]
&\left.-\kappa\lambda_{\hat{n}}\Pi_{\hat{n}}z(t)-\lambda z(t)\right)\mathrm{d}t
+(\sigma(\hat u(t))-\sigma(u(t)))\mathrm{d}W(t),\qquad (t,x)\in \mathbb{R}^{+}\times \mathcal{O},\\[1.0ex]
&\frac{\mathrm{d}}{\mathrm{d}t}\zeta(t)
=\mathcal{T}_\mu\zeta(t)+z(t),\qquad (t,x)\in \mathbb{R}^+\times \mathcal{O},\\[1.0ex]
&(z(0),\zeta(0))
=-U_0\in\mathcal{H}^0_\mu,\qquad x\in\mathcal{O}.
\end{aligned}
\right.
\end{equation}
Define the basic energy functional $\Psi_0(z,\zeta):=\frac12\|z\|_H^2+\frac{1-\kappa}{2}\|\zeta\|_{M_\mu^0}^2.$
Applying It\^{o}'s formula to $\Psi_0(z(t),\zeta(t))$, it yields
\begin{equation}\label{hui}
\begin{aligned}
\mathrm{d}\Psi_0(z(t),\zeta(t))
=&\Bigl[-\kappa a_{\alpha}(z(t),z(t))+(1-\kappa)\langle\mathcal{T}_\mu\zeta(t),\zeta(t)\rangle_{M_\mu^0}\\
&+\langle\varphi(\hat u(t))-\varphi(u(t)),z(t)\rangle_{\langle L^{q}(\mathcal{O}),L^{p_{0}+1}(\mathcal{O})\rangle}\\
&-\kappa\lambda_{\hat n}\|\Pi_{\hat n}z(t)\|_H^2+\frac{1}{2}\|\sigma(\hat u(t))-\sigma(u(t))\|_{\mathcal{L}_{2}(\mathcal{U},H)}^{2}-\lambda\|z(t)\|_{H}^{2}\Bigr]\mathrm{d}t\\[0.8ex]
&+\langle z,\sigma(\hat u(t))-\sigma(u(t))\mathrm{d}W(t)\rangle_{H}.
\end{aligned}
\end{equation}
By the dissipation of $\mathcal{T}_\mu$ in \eqref{lii1}, we arrive at
\[\langle\mathcal{T}_\mu\zeta(t),\zeta(t)\rangle_{M_\mu^0}\leq-\frac{\delta}{2}\|\zeta(t)\|_{M_\mu^0}^2.\]
Moreover, by $(P_3)$, we have
\[\langle\varphi(\hat u(t))-\varphi(u(t)),z(t)\rangle_{\langle L^{q}(\mathcal{O}),L^{p_{0}+1}(\mathcal{O})\rangle}\leq a_\varphi\|z(t)\|_H^2.\]
For the noise term, by assumption $(Q_0)$, we know
\[\|\sigma(\hat u(t))-\sigma(u(t))\|_{\mathcal{L}_{2}(\mathcal{U},H)}^2\leq L_\sigma^{2}\|z(t)\|_H^2.\]
Combining the above estimates, we obtain
\begin{equation}\label{hui1}
\begin{aligned}
&-\kappa a_{\alpha}(z(t),z(t))+(1-\kappa)\langle\mathcal{T}_\mu\zeta(t),\zeta(t)\rangle_{M_\mu^0}+\langle \varphi(\hat u(t))-\varphi(u(t)),z(t)\rangle_{\langle L^{q}(\mathcal{O}),L^{p_{0}+1}(\mathcal{O})\rangle}\\[0.8ex]
&-\kappa\lambda_{\hat n}\|\Pi_{\hat n}z(t)\|_H^2+\frac{1}{2}\|\sigma(\hat u(t))-\sigma(u(t))\|_{\mathcal{L}_{2}(\mathcal{U},H)}^2-\lambda\|z(t)\|_{H}^{2}\\[0.8ex]
\leq&-\kappa a_{\alpha}(\Pi_{\hat n}z(t),\Pi_{\hat n}z(t))+(1-\kappa)\langle\mathcal{T}_\mu\zeta(t),\zeta(t)\rangle_{M_\mu^0}
-\kappa a_{\alpha}((I-\Pi_{\hat n})z(t),(I-\Pi_{\hat n})z(t))\\
&-(\lambda+\kappa\lambda_{\hat n}-a_\varphi-\frac{1}{2}L_\sigma^{2})\|z(t)\|_H^2+\kappa\lambda_{\hat n}\|(I-\Pi_{\hat{n}})z(t)\|_{H}^{2}\\
\leq&-(\lambda+\kappa\lambda_{\hat n}-a_\varphi-\frac{1}{2}L_\sigma^{2})\|z(t)\|_H^2-\frac{\delta(1-\kappa)}{2}\|\zeta(t)\|_{M_\mu^0}^2.
\end{aligned}
\end{equation}
Substituting  \eqref{hui1} to \eqref{hui} and taking  expectation, by $(\Lambda_1)$ and \eqref{shui}, we obtain
\[\frac{\mathrm{d}}{\mathrm{d}t}\mathbb E\Psi_0(z(t),\zeta(t))\leq-\min\left\{2(\lambda+\kappa\lambda_{\hat n}-a_\varphi-\frac{1}{2}L_\sigma^{2}),\delta\right\}\mathbb E\Psi_0(z(t),\zeta(t)).\]
By Gronwall's inequality, we have
\[\mathbb E\Psi_0(z(t),\zeta(t))
\leq e^{-\min\{2(\lambda+\kappa\lambda_{\hat n}-a_\varphi-\frac{1}{2}L_\sigma^{2}),\delta\}t}\Psi_0(z(0),\zeta(0)),\]
which implies
\[\mathbb E\|(z(t),\zeta(t))\|_{\mathcal H_\mu^0}^2
\leq \frac{1}{1-\kappa}e^{-\min\{2(\lambda+\kappa\lambda_{\hat n}-a_\varphi-\frac{1}{2}L_\sigma^{2}),\delta\}t}\|(z(0),\zeta(0))\|_{\mathcal H_\mu^0}^2.\]
This completes the proof of \eqref{chen5}.

We now proceed to prove \eqref{chen6}. Denote by $\widehat{\mathcal L}$ the generator corresponding to \eqref{chen1}. Observe that
\[\widehat{\mathcal L}\Psi_2(\hat u,\hat \eta)=\mathcal L\Psi_2(\hat u,\hat\eta)-\kappa\lambda_{\hat n}\langle \Pi_{\hat n }z,\hat u\rangle_{\mathbb{H}^{2\alpha}}.\]
Applying now Young's inequality,
\begin{align*}
-\kappa\lambda_{\hat n}\langle \Pi_{\hat n}z,\hat u\rangle_{\mathbb{H}^{2\alpha}}&=-\kappa\lambda_{\hat n}\|\Pi_{\hat n}\hat u\|_{\mathbb{H}^{2\alpha}}^2+\kappa\lambda_{\hat n}\langle \Pi_{\hat n}u,\Pi_{\hat n}\hat u+(I-\Pi_{\hat n})\hat u \rangle_{\mathbb{H}^{2\alpha}}\\[0.8ex]
&=-\kappa\lambda_{\hat n}\|\Pi_{\hat n}\hat u\|_{\mathbb{H}^{2\alpha}}^2+\kappa\lambda_{\hat n}\langle \Pi_{\hat n}\hat u, \Pi_{\hat n}u\rangle_{\mathbb{H}^{2\alpha}}\\[0.8ex]
&\leq\frac\kappa4\lambda_{\hat n}\|\Pi_{\hat n}u\|_{\mathbb{H}^{2\alpha}}^2\\[0.8ex]
&\leq\frac\kappa4\lambda_{\hat n}^{3}\|u\|_{H}^2.
\end{align*}
Hence,
\begin{align*}
\widehat{\mathcal L}\Psi_2(\hat u,\hat \eta)\leq \mathcal L\Psi_2(\hat u,\hat\eta)+\frac{\kappa}{4}\lambda_{\hat n}^{3}\|u\|_{H}^2.
\end{align*}
On the one hand, we proceed  exactly as in the proof of Lemma \ref{xiao} to derive a corresponding estimate for $(\hat u(t),\hat\eta(t))$. Namely,
\begin{equation}\label{5.31}
\frac{\mathrm{d}}{\mathrm{d}t}\mathbb{E}\Psi_2((\hat u(t),\hat\eta(t)))\leq-c\mathbb{E}\Psi_2((\hat u(t),\hat\eta(t)))+
C\mathbb{E}\|u(t)\|_H^2+C.
\end{equation}
On the other hand, by Lemma \ref{Z5}, we know that
\[\mathbb{E}\|u(t)\|_H^2\leq 2e^{-\hat{c}_{0}t}\Psi_{0}(u_{0},\eta_0)+2\hat{C}_{0}.\]
%Moreover, for $\mathbb{E}\Psi_{k-1}(\hat U(t))$, by completely repeating the inductive argument in Lemma \ref{xiao}, we infer that
%\[\mathbb{E}\Psi_{k-1}(\hat U(t))\leq C_{k-1}e^{-c_{k-1}t}\Psi_{k-1}(\hat{U}_{0})+C_{k-1}.\]
Since $(\hat u(0),\hat\eta(0))=(0,0)$, substituting these estimates into \eqref{5.31} and  applying Gronwall's inequality, we arrive at
\[\mathbb{E}\Psi_2(\hat u,\hat \eta)\leq C\bigl(\|(u_0,\eta_0)\|_{\mathcal{H}_\mu^{0}}^{2}+1\bigr),
\qquad t\geq 0,\]
which yields \eqref{chen6}. The proof of this lemma is complete.
\end{proof}
Now, we prove  Theorem \ref{chen}.
\begin{proof}
Let $\nu_{\mu}$ be an invariant measure for the Markov semigroup associated with \eqref{lii10}. We prove that
\[\nu_{\mu}\bigl(\mathbb{H}^{2\alpha}\times M_{\mu}^{2\alpha}\bigr)=1.\]
For $R,N\geq1$, define
\[\psi_{R,N}(U):=R\wedge \Bigl(\|\Pi_N u\|_{\mathbb{H}^{2\alpha}}+\|\Pi_N\eta\|_{M_{\mu}^{2\alpha}}\Bigr),\qquad U=(u,\eta)\in \mathcal{H}_{\mu}^{0}.\]
Since $\nu_{\mu}$ is invariant, for every $t\ge 0$ we have
\begin{align}\label{qi10}
\int_{\mathcal{H}_{\mu}^{0}}\psi_{R,N}(U)\nu_{\mu}(\mathrm{d}U)=\int_{\mathcal{H}_{\mu}^{0}} P_t^{\mu}\psi_{R,N}(U)\nu_{\mu}(\mathrm{d}U),
\end{align}
where $P_t^{\mu}$ denotes the Markov transition semigroup associated with the original system  \eqref{lii10}.

Now let $U(t;(u_0,\eta_{0}))=(u(t;u_{0}),\eta(t,\cdot;\eta_{0}))$ be the solution to the original system \eqref{lii10} with initial data $U_0$, and let
$\widehat U(t;0)=(\widehat u(t;0),\widehat \eta(t,\cdot;0))$ be the solution to the controlled system introduced in \eqref{chen1}. Then
\begin{align*}
P_t^{\mu}\psi_{R,N}(U_0)=\mathbb E\psi_{R,N}(U(t;U_0))
\leq&\mathbb E\|\Pi_N(u(t)-\widehat u(t))\|_{\mathbb{H}^{2\alpha}}+\mathbb E\|\Pi_N(\eta(t)-\widehat\eta(t))\|_{M_{\mu}^{2\alpha}}\\[0.8ex]
&+\mathbb E\|\Pi_N\widehat u(t)\|_{\mathbb{H}^{2\alpha}}+\mathbb E\|\Pi_N\widehat\eta(t)\|_{M_{\mu}^{2\alpha}}\\[0.8ex]
\leq&\lambda_{N}\mathbb{E}\|\left(u(t)-\hat{u}(t),\eta(t)-\hat{\eta}(t)\right)\|_{\mathcal{H}_{\mu}^{0}}+\mathbb{E}\|\left(\hat{u}(t),\hat{\eta}(t)\right)\|_{\mathcal{H}_{\mu}^{2\alpha}}.
\end{align*}
By Lemma \ref{chen4}, we have
\[P_t^{\mu}\psi_{R,N}(U_0)\leq C\lambda_{N}e^{-ct}\|U_{0}\|_{\mathcal{H}_{\mu}^{0}}+C\|U_{0}\|_{\mathcal{H}_{\mu}^{0}}^{2}+C.\]
Integrating the above inequality with respect to $\nu_{\mu}$, by \eqref{qi1}, we derive
\begin{equation*}
\begin{split}
\int_{\mathcal{H}_{\mu}^{0}} P_t^{\mu}\psi_{R,N}(U_0)\nu_{\mu}(\mathrm{d}U_0)
&\leq C\lambda_{N}e^{-ct}\int_{\mathcal{H}_{\mu}^{0}}\|U_{0}\|_{\mathcal{H}_{\mu}^{0}}\nu_{\mu}(\mathrm{d}U_0)+C\int_{\mathcal{H}_{\mu}^{0}}\|U_{0}\|_{\mathcal{H}_{\mu}^{0}}^{2}\nu_{\mu}(\mathrm{d}U_0)+C\\[0.8ex]
&\leq C\lambda_{N}e^{-ct}+C.
\end{split}
\end{equation*}
For fixed $N$ and $R$, by choosing $t$ sufficiently large and \eqref{qi10}, we deduce
\[\int_{\mathcal{H}_{\mu}^{0}}\psi_{R,N}(U_0)\nu_{\mu}(\mathrm{d}U_0)\leq C.\]
Furthermore,
\[\int_{\mathcal{H}_{\mu}^{0}}R\wedge \left(\|\Pi_N u\|_{\mathbb{H}^{2\alpha}}+\|\Pi_N\eta\|_{M_{\mu}^{2\alpha}}\right)\nu_{\mu}(\mathrm{d}u,\mathrm{d}\eta)\leq C.\]
Letting $R\to\infty$ and applying the monotone convergence theorem, we infer that
\[\int_{\mathcal{H}_{\mu}^{0}}
\Bigl(\|\Pi_Nu\|_{\mathbb{H}^{2\alpha}}+\|\Pi_N\eta\|_{M_{\mu}^{2\alpha}}\Bigr)\nu_{\mu}(\mathrm{d}u,\mathrm{d}\eta)
\leq C.\]
Furthermore, by the definition of $\Pi_N$ and  applying the monotone convergence theorem once again and letting
$N\to\infty$, we obtain
\[\int_{\mathcal H_\mu^0}\left(\|u\|_{\mathbb{H}^{2\alpha}}+\|\eta\|_{M_\mu^{2\alpha}}\right)\nu^\mu(\mathrm{d}u,\mathrm{d}\eta)\leq C.\]
That is,
\[\nu_{\mu}\bigl(\mathbb{H}^{2\alpha}\times M_{\mu}^{2\alpha}\bigr)=1.\]
This completes the proof.
\end{proof}
\begin{rem}
The restriction to the regularity level $\mathbb{H}^{2\alpha}$ is mainly due to the treatment of the nonlinear term under the present assumptions.
Since
\[\begin{aligned}
\langle \varphi(u),A_\alpha^2u\rangle  &= \langle A_\alpha^{1/2}\varphi(u), A_\alpha^{3/2}u\rangle_H
    \leq\varepsilon\|A_\alpha^{3/2}u\|_H^2+C_\varepsilon\|A_\alpha^{1/2}u\|_H^2\\
    &\leq\varepsilon\|A_\alpha^{3/2}u\|_H^2 +C_\varepsilon\|A_\alpha u\|_H^2.
\end{aligned}
\]
Thus, assumption $(P_4)$ is naturally compatible with the $\mathbb{H}^{2\alpha}$-energy estimate.
In contrast, if one attempts to derive regularity at the level $\mathbb{H}^{3\alpha}$ or higher, one has to control terms involving, for
instance, $A_\alpha\varphi(u)$ and so on. Therefore, under the present structural assumptions,
the $\mathbb{H}^{2\alpha}$ level is the natural regularity level at which the
higher order energy estimate can be closed without imposing further
conditions on the nonlinearity. We also emphasize that $(P_4)$ is not in conflict with
$(P_0)$-$(P_3)$. For example, let
\[\varphi(s)=-as+b\sin s,~a>0,~b\in\mathbb R,\]
where $\varphi$ is globally Lipschitz and $\varphi(0)=0$. Denote by $L_\varphi>0$ a global Lipschitz constant of $\varphi$, i.e.,
\[|\varphi(r)-\varphi(s)|\leq L_\varphi|r-s|,~ r,s\in\mathbb R.\]
then
\[\|A_\alpha^{1/2}\varphi(u)\|_H\leq L_{\varphi}\|A_\alpha^{1/2}u\|_H.\]
Consequently,
\[ s\varphi(s)=-as^2+bs\sin s\leq-\frac{a}{2}s^2+\frac{b^2}{2a},\]
so that $(P_2)$ holds with $p_0=1$. Assumptions $(P_1)$ and $(P_3)$
are also satisfied, while the global Lipschitz continuity of $\varphi$
implies $(P_4)$ by the estimate above. Hence the assumptions
$(P_0)$-$(P_4)$ are mutually compatible.
\end{rem}

\section*{Declaration of competing interest}
\quad\quad There is no competing interest.
\section*{Data availability}
\quad\quad No data was used for the research described in the article.

\end{document}